\documentclass[10pt]{article}
\usepackage{amsmath}
\usepackage{latexsym}
\usepackage{amssymb}
\usepackage{amsfonts}
\usepackage{amsthm}
\title{MACAULAY INVERSE SYSTEMS REVISITED}
\author{J.F. Pommaret \\ CERMICS, Ecole Nationale des Ponts et
  Chauss\'ees,\\6/8 Av. Blaise Pascal, 77455 Marne-la-Vall\'ee Cedex 02,
  France \\
e-mail: pommaret@cermics.enpc.fr }
\date{  }
\begin{document}
\maketitle

\noindent
{\bf ABSTRACT}:\\

         Since its original publication in 1916 under the title {\it The Algebraic Theory of Modular Systems}, the book [13] by F.S. Macaulay has attracted a lot of scientists with a view towards pure matematics [6] or applications to control theory [15] through the last chapter dealing with the so-called {\it inverse system}. The basic intuitive idea is the well known parallel existing between ideals in polynomial rings and systems of partial differential (PD) equations in one unknown with constant coefficients. \\
       A first purpose of this paper is thus to extend these results to arbitrary systems of PD equations by exhibiting a link with the formal theory of systems of PD equations ([17],[23],[24]) where concepts such as {\it involution} are superseding the H-{\it bases} of Macaulay.\\
       The second idea is to transfer the properties of ideals to their residue modules, in particular to extend to differential modules the {\it unmixedness} assumption of Macaulay. For this we use extensively the results of modern {\it algebraic analysis} ([2],[11],[16],[18],[19]), revisiting in particular the concept of {\it purity} by means of {\it localization} techniques. Accordingly, this paper can also be considered as a refinement and natural continuation of [20].\\
       Finally, following again Macaulay in the differential setting, the cornerstone and main novelty of the paper is to replace the {\it socle} of a module by the {\it top} of the corresponding dual system in order to be able to look for generators by using well known arguments of algebraic geometry such as Nakayama's lemma [12].\\
        Many explicit examples are provided in order to illustrate the main constructive results that provide new hints for applying computer algebra to algebraic analysis [21].\\
      

\noindent
1) {\bf INTRODUCTION}:\\

  With only a slight abuse of language, one can say that the birth of the {\it formal theory} of systems of ordinary differential (OD) equations or partial differential (PD) equations is coming from the work of M. Janet in 1920 [9] along algebraic ideas brought by D. Hilbert at the same time in his study of sygyzies for finitely generated modules over polynomial rings. In 1965 [5] B. Buchberger invented {\it Gr\"{o}bner bases}, named in honor of his Phd advisor W. Gr\"{o}bner, whose earlier work done in 1940 on polynomial ideals and PD equations with constant coefficients provided a source of inspiration [8]. However, the approaches of Janet and Buchberger/Gr\"{o}bner both suffer from the same lack of {\it intrinsicness} as they highly depend on the ordering of the $n$ independent variables and derivatives of the $m$ unknowns [7,23]. \\
   
     Meanwhile, {\it commutative algebra}, namely the study of modules over rings, was facing a very subtle problem, the resolution of which led to the modern but difficult {\it homological algebra} with sequences and diagrams. Roughly, one can say that the problem was essentially to study properties of finitely generated modules not depending on the {\it presentation} of these modules by means of generators and relations. This very hard step is based on homological/cohomological methods like the so-called {\it extension} modules which cannot therefore be avoided ([4],[18],[22]).\\
     
     In order to sketch this problem, let us present two simple examples. We shall use standard notations similar to the ones of computer algebra, namely a dot represents the derivative with respect to a single independent variable (time in classical control theory) while, in the case of many independent variables $(x^1,..., x^n)$, the notation $d_{ij}=d_{ji}$ describes for example the second order derivative with respect to $x^i,x^j$ with $1\leq i,j\leq n$. In the first case with standard notations, everybody will understand at once that integrating the second order OD equation $\ddot{y}=0$ with $m=n=1$ is equivalent to integrating the system of two first order OD equations ${\dot{y}}^1-y^2=0, {\dot{y}}^2=0$. However, even with $m=n=2$ and the same two unknowns $u,v$ in both cases, it is not evident at all that integrating the second order PD equation $d_{12}u-d_{22}v-u=0$ is equivalent to integrating the system of two fourth order PD equations $d_{1122}u-d_{1222}v-d_{22}v-u=0, d_{1112}u-d_{1122}v-d_{11}u=0$.\\
     
     As before, using now rings of {\it differential operators} instead of polynomial rings led to {\it differential modules} and to the challenge of adding the word {\it differential} in front of concepts of commutative algebra. Accordingly, not only one needs properties not depending on the presentation as we just explained but also properties not depending on the coordinate system as it becomes clear from any application to mathematical or engineering physics where tensors and exterior forms are always to be met like in the space-time formulation of electromagnetism. Unhappily, no one of the previous techniques for OD or PD equations could work.\\

     By chance, the intrinsic study of systems of OD or PD equations has been pioneered in a totally independent way by D. C. Spencer and collaborators after 1960 [24], in order to relate differential properties of the PD equations to algebraic properties of their symbols, a technique superseding the {\it leading term} approach of Janet or Gr\"{o}bner .\\
     
     Accordingly, it was another challenge to unify the {\it purely differential} approach of Spencer with the {\it purely algebraic} approach of commutative algebra, having in mind the necessity to use the previous homological algebraic results in this new framework. This sophisticated mixture of differential geometry and homological algebra, now called {\it algebraic analysis}, has been achieved after 1970 by V. P. Palamodov for the constant coefficients case [16], then by M. Kashiwara [11] for the variable coefficients case.\\
     
     It is only in 1990, thanks to the work of U. Oberst, that such a theory has been applied with success to control theory [15]. Then the things went on rather fast towards computer algebra and many packages now exist for computing the extension modules and related concepts ([21] is sufficient for dealing with most of this paper). Of course, many difficult problems are left and we provide details about a few of them, having in mind the recent workshops on Gr\"{o}bner bases and applications successively held at RISC/Linz in 2006 and 2008 [20].\\
     
     When a given system of linear PD equations of order $q$ is given, it defines by residue a differential module $M$ over the underlying ring $D$ of differential operators. Then it becomes today possible to decide by means of computer algebra the class to which $M$ belongs among $n$ classes ranging from {\it free, torsion-free, reflexive}, ... , to {\it projective} and {\it free} [21]. However, the set of elements of $M$, namely the finite linear combinations of the unknowns and their derivatives modulo the given PD equations and their derivatives, such that {\it each of them} does satisfy {\it at least one} PD equation {\it for itself}, provides the {\it torsion submodule} $t(M)$ and $M$ is torsion-free if its torsion submodule is zero. An open but useful problem, independent of the previous classification, is now to classify elements in $t(M)$. For this, we recall that the Hilbert-Serre theorem asserts that the {\it dimension} $d(M)$ of a module $M$ defined by a system of PD equations is equal to the dimension $d(V)$ of the {\it characteristic variety} $V$ of the system and this number does not depend on the presentation and filtration of $M$ ([14],[18], p 542,544). Then we define $t_r(M)$ to be zero or the unique greatest differential submodule of $M$ having dimension $<n-r$ and we have the nested chain of $n$ differential submodules:\\
     \[   0=t_n(M)\subseteq  t_{n-1}(M)\subseteq . . . \subseteq t_1(M)\subseteq t_0(M)=t(M)\subseteq M   \]
A basic question is thus to determine the {\it classes} and the {\it gaps} in the above chain, as indeed, in many known explicit situations, a few intermediate modules do coincide. The interest is to provide new domains of applications and we sketch the underlying idea on a simple academic example.\\

     With the above notations notations and $m=1,n=q=2$, the system $d_{22}y=0, d_{12}y=0$ clearly determines a torsion module $t(M)=M$. The torsion elements $z'=d_1y$ and $z"=d_2y$ do not satisfy the same decoupling PD equations as $z'$ only satisfies $d_2z'=0$ while $z"$ satisfies $d_2z"=0, d_1z"=0$ and we have the nested chain with strict inclusions:\\
     \[   0=t_2(M)\subset t_1(M) \subset t_0(M)=t(M)=M  \]
the classification being obtained through the {\it dimension} $d$ or rather {\it codimension} $cd=n-d$ of the differential modules generated by the respective torsion elements as we have indeed $cd(Dz")=2$ and $cd(Dz')=1$. Of course, the same decoupling type problem can be asked for any engineering system in gasdynamics or magnetohydrodynamics (MHD) but we do not know a single work existing towards such a classification in view of the difficulty of the mathematical framework involved. As an ultimate goal, a particularly important problem should be to study the dependence of the previous classification on parameters when the system depends on certain constant parameters, a result generalizing the controllability problem for OD systems in control theory where $n=1$ only [20].\\
 
    Accordingly, the hope should be to have a computer algebra package providing the classes, the gaps and eventually generating elements. A particularly interesting case should be to characterize $r$-{\it pure} modules, namely modules $M$ such that there exists an integer $0\leq r\leq n$ with $t_r(M)=0$ and $ t_{r-1}(M)=M$. Equivalently, $M$ is $r$-pure whenever $cd(Dm)=r, \forall m\in M$. For constant coefficients systems in one unknown, such a concept had been discovered in 1916 by F. S. Macauly under the name  {\it unmixed ideal} ([13], glossary of the last edition and \S 41,77).\\   
    In fact, while looking at the last chapter of his book since many years, we were convinced that the double picture of p. 67 was nothing else than, {\it sise by side}, the (lower triangular) matrix of the coefficients of the system of OD/PD equations, organized horizontally along the increasing order of the derivatives of the unknowns and vertically along the increasing order of the leading terms of the equations with respect to the previous ordering, combined with an ordering of the various possible formal solutions made up by {\it truncated power series}, the underlying idea being to cancel successively the terms of order zero, then zero and one, ... and so on. However, it is only a few months ago that we suddenly understood the true reason for supposing, {\it as a crucial assumption indeed though it is only presented as a purely technical argument} (p. 89), that the ideal under study was {\it unmixed}. We explain thereafter this point.\\
    First of all, the properties (prime, primary, unmixed,...) attributed to an ideal $\mathfrak{a}$ in the ring $A=k[\chi]=k[{\chi}_1,...,{\chi}_n]$ of polynomials in the $n$ indeterminates ${\chi}_1,...,{\chi}_n$ with coefficients in the field $k$ are now, along with the modern setting of commutative algebra, attributed to the residual module $A/{\mathfrak{a}}$.  Then we got in mind that, in the study of an $r$-pure differential modules, a delicate though expected theorem is stating that the corresponding characteristic variety is {\it unmixed} too, with the same meaning as above, that is the underlying algebraic set is the union of irreducible components of the {\it same codimension} $r$ ([2], p 42,[18], p 551). This analogy was thus giving rise to the challenge of relating the work of Macaulay on unmixed polynomial ideals to the study of purity for differential modules. In particular, the extension to this new framework of a {\it localization criterion}, provided by Macaulay in the classical setting, constitutes one of the main results of this paper and provides new hints for applying computer algebra to algebraic analysis.\\
    
In section 2 we start presenting a few motivations from commutative algebra, then apply localization theory to systems of OD equations and finally generalize the results so far obtained to systems of PD equations.\\

Then section 3 establishes a way to use a partial localization in order to test the purity of a module as a basic assumption replacing the unmixedness of an ideal in the work of Macaulay.\\

The key section 4 describes the inverse system and exhibits the duality existing between the socle of a module and the top of the corresponding system in the sense of Spencer or Macaulay.\\

The final section 5 explains Macaulay's secrete as a way to use Nakayama's lemma in order to find out generating sections of the systems corresponding to pure modules.\\

We end the present section explaining this point on a few simple but illuminating examples. Using a sub-index $x$ for the derivatives when $n=1$, the general solution of $y_{xx}-y=0$ is $y=a e^x+b e^{-x}$ with $a,b$ constants and the derivative of $e^x$ is $e^x$ while the derivative of $e^{-x}$ is $-e^{-x}$. Hence we could believe that we need a basis $\{e^x,e^{-x}\}$ with {\it two} generators for obtaining all the solutions through derivatives. However, setting as usual $sh(x)=\frac{1}{2}(e^x-e^{-x}), ch(x)=\frac{1}{2}(e^x+e^{-x})$, we have equivalently $y=c\times sh(x)+d\times ch(x)$ with $c,d$ constants. As the derivative of $ch(x)$ is $sh(x)$, we need only a basis $\{ch(x)\}$ with {\it one} generator . If we now consider the system $y^1_{xx}=0, y^2_x=0$, we need a basis $\{(x,0),(0,1)\}$ with {\it two} generators. However, changing slightly the latter system to $y^1_{xx}-y^1=0, y^2_x=0$ and introducing $z=y^1-y^2$, it is equivalent to set $y^1=z_{xx}, y^2= z_{xx}-z$ and consider the system $z_{xxx}-z_x=0$ with the only generator $\{ch(x)-1\}$ leading therefore to the only generator $\{(ch(x),1)\}$ for the original system .\\

 \noindent
2) {\bf MOTIVATIONS}\\

Let $k$ be a field of characteristic zero and $\chi=({\chi}_1,...,{\chi}_n)$ be indeterminates over $k$. We introduce the ring $A=k[{\chi}_1,...,{\chi}_n]$ of polynomials with coefficients in $k$ and various classes of ideals. The set of {\it maximum ideals} is denoted by $max(A)$ with elements ${\mathfrak{m}}$,..., the set of (proper) {\it prime ideals} is denoted as usual by $spec(A)$ with elements ${\mathfrak{p}}$,... and the set of {\it primary ideals} with elements ${\mathfrak{q}}$ such that $ab\in \mathfrak{q}, b\notin \mathfrak{q}\Rightarrow a\in \mathfrak{p}=rad(\mathfrak{q})$ that is $a^r\in \mathfrak{q}$ for a certain integer $r\in \mathbb{N}$. The importance of primary ideals lies in the fact, largely emphasized by Macaulay, that any ideal $\mathfrak{a}$ can be written as a finite irredundant intersection $\mathfrak{a}={\mathfrak{q}}_1\cap ... \cap{\mathfrak{q}}_s$ of primary ideals, called {\it primary decomposition}. Setting ${\mathfrak{p}}_i=rad({\mathfrak{q}}_i)$, we obtain at once the {\it prime decomposition} $rad(\mathfrak{a})={\mathfrak{p}}_1\cap ... \cap {\mathfrak{p}}_s$ though sometimes this new decomposition may not be irredundant with strict inclusion ${\mathfrak{p}}_i\subset{\mathfrak{p}}_j$ for certain couples of indices $(i,j)$. In this case one uses to say that the component defined by ${\mathfrak{p}}_j$ is {\it embedded} into the component defined by ${\mathfrak{p}}_i$ in the algebraic set defined by $\mathfrak{a}$. Also, for any prime ideal $\mathfrak{p}$, we denote by $cd(A/ \mathfrak{p})=n-d(A/ \mathfrak{p})$ the {\it codimension} of $A/ \mathfrak{p}$ with $d(A/ \mathfrak{p})=trd(Q(A/ \mathfrak{p})/k)$ the {\it transcendence degree} of the algebraic extension of the field of fractions of the integral domain $A/\mathfrak{p}$ over the field $k$. For an arbitrary ideal $\mathfrak{a}$, the codimension is usually denoting the minimum among the codimensions of the components defined by the minimum prime ideals in the corresponding prime decomposition, which are therefore not embedded.\\

\noindent
{\bf DEFINITION 2.1}: An ideal $\mathfrak{a}\subset A$ is {\it unmixed} if $cd(A/ {\mathfrak{p}}_1)=...=cd(A/ {\mathfrak{p}}_s)$ in a primary decomposition and we have therefore ${\mathfrak{p}}_i \nsubseteq{\mathfrak{p}}_j, \forall (i,j)$. Otherwise $\mathfrak{a}$ is said to be {\it mixed}.\\

We now present a few examples that will be used in the sequel with a totally different approach.\\

\noindent
{\bf EXAMPLE 2.2}: $\mathfrak{q}=(({\chi}_3)^2,{\chi}_1{\chi}_3-{\chi}_2)$ is primary with $rad(\mathfrak{q})=\mathfrak{p}=({\chi}_3,{\chi}_2)\Rightarrow cd(A/ \mathfrak{q})=2$. Similarly $\mathfrak{a}=({\chi}_1,{\chi}_2{\chi}_3)=({\chi}_1,{\chi}_2)\cap ({\chi}_1,{\chi}_3) $ is unmixed with $cd(A/ \mathfrak{a})=2$ but the new ideal $\mathfrak{a}=(({\chi}_1)^2,{\chi}_1{\chi}_2,{\chi}_1{\chi}_3,{\chi}_2{\chi}_3)=({\chi}_1,{\chi}_2)\cap ({\chi}_1,{\chi}_3)\cap ({\chi}_1,{\chi}_2,{\chi}_3)^2$ is mixed with two minimum prime ideals and one embedded component. More generally, any ideal having a basis containing as many polynomials as the codimension of the corresponding residual module has been called {\it ideal of the principal class} by Macaulay who proved that any such ideal is unmixed ([13], \S 48, p 40,49). For a modern approach through regular sequences, see ([12] , VI, 3, p 183).\\

\noindent
{\bf EXAMPLE 2.3}: $\mathfrak{a}=(({\chi}_2)^2,{\chi}_1{\chi}_2)=({\chi}_2)\cap ({\chi}_1,{\chi}_2)^2 ={\mathfrak{q}}_1\cap{\mathfrak{q}}_2$ is mixed with ${\mathfrak{q}}_1={\mathfrak{p}}_1\subset {\mathfrak{p}}_2=rad({\mathfrak{q}}_2)=({\chi}_1,{\chi}_2)\in max({\mathbb{Q}}[{\chi}_1,{\chi}_2])$.\\

\noindent
{\bf EXAMPLE 2.4}: ([13], \S 42, p 44)  $\mathfrak{a}=(({\chi}_1)^3,({\chi}_2)^3,(({\chi}_1)^2+({\chi}_2)^2){\chi}_4+{\chi}_1{\chi}_2{\chi}_3)$ is mixed with $s=4$ and only one minimum prime because $({\chi}_1{\chi}_2)^2({\chi}_1,{\chi}_2,{\chi}_3,{\chi}_4)\subset {\mathfrak{a}} $ but ${\chi}_1{\chi}_2\notin \mathfrak{a}$.\\

\noindent
The main idea is then to {\it transfer} the properties of an ideal $\mathfrak{a}\subset A$ to the residue module $M=A/\mathfrak{a}$ over $A$.\\

\noindent
{\bf DEFINITION 2.5}: A module $M$ is said to be {\it prime} ({\it primary}) if $ax=0, 0\neq x\in M\Rightarrow aM=0$ ($a^rM=0$ for a certain integer $r$) though people sometimes add the prefix "co". \\

Now, having in mind the so-called {\it chinese remainder theorem} ([12], p 41), any primary decomposition gives rise to a monomorphism $0\rightarrow M \rightarrow Q_1\oplus ... \oplus Q_s$ with primary modules $Q_i=A/{\mathfrak{q}}_i$ and epimorphisms $M\rightarrow Q_i\rightarrow 0, \forall i=1,...,s$. Conversely, looking for such a situation allows to exhibit a primary decomposition for {\it reducible} modules (see [18], p 110 for more details).\\

It is now tempting, and this too was a key idea of Macaulay, to introduce $n$ commuting derivatives $d_1,...,d_n$ for which $k$ should be a field of constants and to introduce the ring $D=k[d]=k[d_1,...,d_n]$ of differential operators with coefficients in $k$. As $D$ and $A$ are isomorphic by $d_i\leftrightarrow {\chi}_i$, any ({\it nonlinear}) ideal of $A$ gives rise to a ({\it linear}) system of OD/PD equations {\it in one unknown only} and conversely. {\it It thus remains to use techniques for} PD {\it equations in order to study ideals or modules}. However, the situation for a differential field $K$ with subfield of constants $k$ and/or systems of PD equations for many unknowns escapes from the previous approach and we conjectured that they could be treated {\it by their own}, the specific situation considered by Macaulay being just a particular case of the general theory that we shall present in this paper.\\

First of all, we sketch the technique of {\it localization} in the case of OD equations, comparing to the situation met in classical control theory where $n=1$ and the dimension can therefore only be 0 or 1. For this, setting as usual $d=d_1=d/dt=dot$, we may introduce (formal) unknowns $y^1,...,y^m$ and set $Dy=Dy^1+...+Dy^m\simeq D^m$. If we have a given system $\Phi=0$ of OD equations of order $q$, a basic question in control theory is to decide whether the control system is "{\it controllable}" or not. It is not our purpose to discuss here about such a question (see [10],[15],[18],[19] for more details) but we just want to state the final formal test in terms of a property of the differential module $M=Dy/D\Phi$. Care must be taken that in the sequel, for simplicity and unless needed, we shall not always put a "{\it bar}" on the residual image of $y$ in the canonical projection $Dy\rightarrow M$. We explain our goal on an example.\\

\noindent
{\bf EXAMPLE 2.6}: With $m=3$ and a constant parameter $a$, we consider the first order system ${\Phi}^1\equiv {\dot{y}}^1-ay^2-{\dot{y}}^3=0, {\Phi}^2\equiv y^1-{\dot{y}}^2+{\dot{y}}^3=0$. Let us apply Laplace transform $\hat{y}(s)={\int}^{\infty}_0e^{st}y(t)dt$ to this system. Using the integration by part formula ${\int}^{\infty}_0e^{st}{\dot{y}}(t)dt=[e^{st}y(t)]^{\infty}_0-s\hat{y}(s)$, we should eventually need to know $y(0)$ though the Kalman test of controllability is purely formal as it only deals with ranks of matrices [10]. Since a long time we had in mind that setting $y(0)=0$ was not the right way and that Laplace transform could be superseded by another purely formal technique. For this, let us replace {\it formally} $d$ by the purely algebraic symbol $\chi$ whenever it appears and obtain the system of {\it linear equations} :\\
\[ \chi y^1-ay^2-\chi y^3=0, 1y^1-\chi y^2+\chi y^3=0 \Rightarrow y^1=\frac{\chi (\chi +a)}{{\chi}^2-a} y^3, y^2=\frac{\chi (\chi +1)}{{\chi}^2-a} y^3 \]
but we could have adopted a different choice for the only arbitrary unknown. At this step there are only two possibilities :\\
$\bullet a\neq 0,1 \Rightarrow $no {\it simplification} may occur and, getting rid of the common denominator, we obtain an algebraic parametrization leading to a differential parametrization as follows:\\
\[ y^1=\chi (\chi +a) z, y^2=\chi (\chi +1) z, y^3=({\chi}^2-a)z \Rightarrow y^1=\ddot{z}+a\dot{z}, y^2=\ddot{z}+\dot{z}, y^3=\ddot{z}-az  \]
$\bullet a=0$ or $a=1 \Rightarrow$  a {\it simplification} may occur. For example, with $a=0$, setting $z=y^1-y^3$ we obtain $\chi z=0$ that is to say $\dot{z}=0$.\\

Recapitulating, we discover that a control system is controllable iff one cannot get any {\it autonomous} element satisfying an OD equation {\it by itself}. For understanding such a result in an algebraic manner, let $M$ be a module over an integral domain $A$ containing 1. A subset $S\subset A$ is called a {\it multiplicative subset} if $1\in S$ and $\forall s,t \in S \Rightarrow st \in S$. Moreover, we shall need and thus assume the {\it Ore condition} on $S$ and $A$, namely $aS\cap sA\neq \emptyset , 
\forall a\in A, s\in S$.\\

\noindent
{\bf DEFINITION 2.7}: For any module $M$ over $A$, we define $S^{-1}M=\{s^{-1}x{\mid}s\in S,x\in M/\sim\}$ with $s^{-1}x\sim t^{-1}y \Leftrightarrow \exists u,v\in A, us=vt\in S , ux=vy$. We have $S^{-1}M=S^{-1}A{\otimes}_AM$ and we set $t_S(M)=\{x\in M\mid \exists s\in S, sx=0\} $ in the exact sequence $0\rightarrow t_S(M)\rightarrow M \rightarrow S^{-1}M$ where  the last morphism is $x \rightarrow 1^{-1}x$.\\

\noindent
{\bf EXAMPLE 2.8}: $S=A-\{ 0 \} \Rightarrow S^{-1}A=Q(A)=K $ field of fractions of $A$ and we introduce the {\it torsion submodule} $t_S(M)=t_0(M)=t(M)=\{ x\in M \mid \exists 0\neq a\in A, ax=0\}$ of $M$. Also, if $\mathfrak{p}\in spec(A)$ and $S=A-\mathfrak{p}$, one uses to set $S^{-1}M=M_{\mathfrak{p}}$.\\

\noindent
{\bf PROPOSITION 2.9}: When $M$ is finitely generated and $t(M)=0$, from the inclusion $M\subset K{\otimes}_A M$, we deduce that there exists a finitely generated free module $F$ with $M\subset F$.\\ 

\noindent
{\bf REMARK 2.10}: Though the above proposition provides a parametrization for any $n$ in the case of a torsion-free module, in the particular case $n=1$ there is an isomorphism $M\simeq t(M)\oplus M/t(M)$ not so well known in OD control theory. Indeed the projection onto the second factor is the canonical projection onto the torsion-free module $M/t(M)$ which is a free and thus projective module when $n=1$, a result allowing to split the short exact sequence $0\rightarrow t(M)\rightarrow M \rightarrow M/t(M) \rightarrow 0$. This is not evident at all on Example 2.6 and even on the very simple example ${\dot{y}}^1-{\dot{y}}^2=0$.\\

The comparison with Example 2.6 needs no comment at least when $n=1$ and controllability must have to do with $t(M)=0$ when $n\geq 2$ though it is only quite later on in the paper that we shall be able to generalize the result expressed by the above remark. Also the extension of the above results to the non-commutative case $D=K[d]$ where $K$ is a differential field with $n$ commuting derivations ${\partial}_1,...,{\partial}_n$  can be achieved but is much more delicate ([11],[18],[19]).\\

\noindent
{\bf EXAMPLE 2.11}: When $a=a(t)$ in Example 2.6, the controllability condition is now the Ricatti inequality $\dot{a} + a^2-a\neq 0$ in a coherent way with the constant coefficient case already considered.\\

Taking into account the works of Janet and Spencer, the study of systems of PD equations cannot be achieved without understanding {\it involution} and we now explain this concept (compare to [7,23]). For this, let $\mu=({\mu}_1,...,{\mu}_n)$ be a multi-index with {\it length} $\mid\mu\mid={\mu}_1+ ... +{\mu}_n$. We set $\mu +1_i=({\mu}_1,...,{\mu}_{i-1},{\mu}_i+1,{\mu}_{i+1},...,{\mu}_n)$ and we say that $\mu$ is of {\it class} $i$ if ${\mu}_1=...={\mu}_{i-1}=0, {\mu}_i\neq 0$. Accordingly, any operator $P=a^{\mu}d_{\mu}\in D$ acts on the unknowns $y^k$ for $k=1,...,m$ as we may set $d_{\mu}y^k=y^k_{\mu}$ with $y^k_0=y^k$ and introduce the {\it jet coordinates} $y_q=\{y^k_{\mu}\mid k=1,...,m; 0\leq \mid \mu \mid \leq q\}$. It follows that, if a system of PD equations can be written in the form ${\Phi}^{\tau}\equiv a^{\tau\mu}_ky^k_{\mu}=0$ with $a\in K $, we may introduce the differential module $M=Dy/D\Phi$ but we notice that the work of Macaulay only covers the case $m=1$. Then we define the (formal) {\it prolongation} of ${\Phi}^{\tau}$ with respect to $d_i$ to be $d_i{\Phi}^{\tau}\equiv a^{\tau\mu}_ky^k_{\mu +1_i}+ {\partial}_ia^{\tau\mu}_ky^k_{\mu}$ and induce maps $d_i:M\rightarrow M: {\bar{y}}^k_{\mu}\rightarrow {\bar{y}}^k_{\mu +1_i}$ by residue. It follows that the module $M$ is endowed with a quotient filtration induced by its presentation which is a {\it strict} morphism when the defining system is involutive ([2],[14], p 383,[18], p 445).\\

Changing linearly the derivations if necessary, we may successively solve the maximum number of equations with respect to the jets of class $n$, class $(n-1)$,..., class 1. At each order, a certain number of jets called {\it principal} ($pri$) can therefore be expressed by means of the other jets called {\it parametric} ($par$). Moreover, for each equation of order $q$ and class i, $d_1,...,d_i$ are called {\it multiplicative} while $d_{i+1},...,d_n$ are called {\it nonmultiplicative} and $d_1,...,d_n$ are nonmultiplicative for all the remaining equations of order $\leq q-1$ (Pommaret basis in [23]).\\

\noindent
{\bf DEFINITION 2.12}: The system is said to be {\it involutive} if each prolongation with respect to a nonmultiplicative derivation is a linear combination of prolongations with respect to the multiplicative ones. Using Spencer cohomology, one can prove that such a definition is in fact intrinsic [7,17,18,23,24] and generalizes the concept of H-{\it basis} used by Macaulay ([13], p 36,39,67,68,86).\\

\noindent
{\bf EXAMPLE 2.13}: ([13], \S 38, p 40 where one can find the first intuition of formal integrability) The ideal $\mathfrak{q}=(({\chi}_1)^2, {\chi}_1{\chi}_3-{\chi}_2)$ provides the system $y_{11}=0, y_{13}-y_2=0$ which is not involutive. Effecting the permutation $(1,2,3)\rightarrow (3,2,1)$, we get the new system $y_{33}=0, y_{13}-y_2=0$. As $d_1y_{33}-d_3(y_{13}-y_2)=y_{23}$ and $d_1y_{23}-d_2(y_{13}-y_2)=y_{22}$, the new system $y_{33}=0, y_{23}=0, y_{22}=0, y_{13}-y_2=0$ is involutive with 1 equation of class 3, 2 equations of class 2 and 1 equation of class 1.\\

\noindent
{\bf APPLICATION 2.14}: $t(M)=M$ iff the number of equations of class $n$ is $m$. Otherwise there is a strict inclusion $t(M)\subset M$.\\

\noindent
{\bf PROPOSITION 2.15}: ([17,24]) The following recipe will allow to bring an involutive system of order $q$ to an equivalent (isomorphic modules) involutive system of order 1 {\it with no zero order equations} called 
{\it Spencer form}:\\
1) Use all parametric jets up to order $q$ as new unknowns.\\
2) Make one prolongation.\\
3) Substitute the new unknowns.\\

\noindent
{\bf PROPOSITION 2.16}: For such a system, defining the {\it character} ${\alpha}^i_1=m-$number of equations of class i, we have ${\alpha}^1_1\geq {\alpha}^2_1\geq ... \geq {\alpha}^n_1\geq 0$. The first nonzero character and the number of nonzero characters are intrinsic integers, coming from the Hilbert polynomial of the module/system (See [23] for a computer algebra implementation).\\

\noindent
{\bf REMARK 2.17}: $cd_D(M)\geq r \Leftrightarrow {\alpha}^{n-r}_1\neq 0,{\alpha}^{n-r+1}_1=...={\alpha}^n_1=0$. In that case, we shall say that we have {\it full class} n,..., {\it full class} (n-r+1). Thus $r=2$ in the above example.\\

\noindent
3) {\bf PURITY}\\

As a first basic fact, quite important for the study of the noncommutative case, one must {\it carefully} distinguish between an ideal/system and its {\it symbol} part, namely the top order part of order $q$ when the system is involutive. The following example will explain the difficulty involved, hidden in the use of Gr\"{o}bner bases which are not intrinsically defined.\\

\noindent
{\bf EXAMPLE 3.1}: The primary ideal corresponding to the involutive system of the preceding example of Macaulay has radical $\mathfrak{p}=({\chi}_3,{\chi}_2)$. However, the symbol part is the homogeneous ideal $\mathfrak{a}=(({\chi}_3)^2, {\chi}_2{\chi}_3, ({\chi}_2)^2, {\chi}_1{\chi}_3)=(({\chi}_2)^2, {\chi}_3)\cap ({\chi}_1, {\chi}_2, {\chi}_3)^2$ with the same radical but it is a pure coincidence. Nevertheless, $\mathfrak{a}$ and $\mathfrak{q}$ have the same (co)dimension according to the famous Hilbert-Serre theorem that we shall recall in Proposition 3.6. The importance of involution has not been pointed out clearly in the study of Gr\"{o}bner bases [20]. However, it is clear that the two polynomials $({\chi}_3)^2$ and ${\chi}_1{\chi}_3-{\chi}_2$ generate $\mathfrak{q}$ in the previous example but the corresponding homogeneous ideal at order 2 is now $(({\chi}_3)^2, {\chi}_1{\chi}_3)$ with radical $({\chi}_3)$ providing a wrong dimension and the corresponding presentation is no longer strict.\\

\noindent
{\bf DEFINITION 3.2}: The {\it characteristic variety} $V$ of an involutive system of order $q$ is the algebraic set defined by the {\it radical} (care) of the polynomial ideal in $K[\chi]$ generated by the $m\times m$ minors of the {\it characteristic matrix} $(a^{\tau\mu}_k{\chi}_{\mu})$ where $\mid\mu\mid =q$. Of course, when $m=1$ we recover the radical of the symbol ideal of the ideal we started with and the involutive assumption is essential.\\

If $m\in M$, then the differential submodule $Dm\subset M$ is defined by a system of OD/PD equations for one unknown only and we may look for its codimension $cd_D(Dm)$. In the commutative case, looking at the annihilators, we get $ann(M)\subset ann(Dm)$. In particular, if $M$ is primary its annihilator is a primary ideal $\mathfrak{q}$ with radical $\mathfrak{p}$ and we have $\mathfrak{q} \subseteq ann(Dm) \subseteq \mathfrak{p}, \forall m\in M$ as a possible characterisation. Accordingly, if $M$ is prime, then $ann(Dm)=\mathfrak{p}, \forall m\in M$.\\

\noindent
{\bf EXAMPLE 3.3}: In Example 2.13, with the primary ideal $\mathfrak{q}$, then $y_2$ and $y_3$ are killed by $\mathfrak{p}$ though $y$ is killed by $\mathfrak{q}$. The situation changes completely with the corresponding homogeneous ideal $\mathfrak{a}$ as $y_1$ is killed by $({\chi}_3,({\chi}_2)^2)$ and $y_{12}$ is killed by $({\chi}_3, {\chi}_2)$ though $y_3$ is killed by $({\chi}_1,{\chi}_2,{\chi}_3)$.\\

Even in the noncommutative case of systems with coefficients in a differential field $K$, one can prove with the homological techniques of algebraic analysis using extension modules and bidualizing complexes ([2], Theorem A.IV.2.14, p 494,[18], Proposition IV.3.161, p 545) :\\

\noindent
{\bf PROPOSITION 3.4}: $t_r(M)=\{m\in M\mid cd(Dm)>r\}$ is the greatest differential submodule of $M$ having codimension $>r$.\\

\noindent
{\bf PROPOSITION 3.5}: $t_r(M)$ does not depend on the presentation and thus on the filtration of the module $M$ as it can be defined inductively by the exact sequences :\\
\[  0 \longrightarrow  t_r(M) \longrightarrow t_{r-1}(M) \longrightarrow ext^r_D(ext^r_D(M,D),D)  \]
if we start from $t_{-1}(M)=M$ and $t_0(M)=t(M)$ when $r=0$.\\

Thanks to its implementation in [21], this proposition is essential for the use of computer algebra and allows to refer to the Spencer form. In fact, {\it the situation is exactly similar to that of control theory with the Kalman form}.\\

\noindent
{\bf PROPOSITION 3.6}: $cd_D(M)=cd(V)=r \Leftrightarrow {\alpha}^{n-r}_1\neq 0, {\alpha}^{n-r+1}_1= ... ={\alpha}^n_1=0 \Leftrightarrow t_r(M)\neq M,t_{r-1}(M)= ... =t_0(M)=t(M)=M$. \\

We may therefore define as in the Introduction:\\

\noindent
{\bf DEFINITION 3.7}: $M$ is r-{\it pure} $\Leftrightarrow t_r(M)=0, t_{r-1}(M)=M, \Leftrightarrow cd(Dm)=r, \forall m\in M$. In particular, $M$ is 0-pure iff $t(M)=0$. Otherwise, if $cd_D(M)=r$ but $M$ is not pure, by analogy with Remark 2.10, we shall call $M/t_r(M)$ the {\it pure part} of $M$.\\

The following key results using a kind of {\it partial localization} generalize the similar ones first obtained by Macaulay ([13], \S 82) and provide a technical test linking purity and involution, both with an effective construction of Proposition 3.4. From now on we shall only consider the constant coefficient situation, considering ${\chi}_1, ..., {\chi}_{n-r}$ just like parameters ([13], \S 77, p 86), but most of the results can be extended to the variable coefficient situation, though with a lot of work more.\\

\noindent
{\bf THEOREM 3.8}: If $cd_D(M)=r$ one has the exact sequence:\\
\[    0 \longrightarrow t_r(M) \longrightarrow M \longrightarrow k({\chi}_1,...,{\chi}_{n-r})\otimes M  \]
\noindent
{\it Proof}: According to the definition of involution, the system made by the PD equations of class 1+...+class $(n-r)$ is also involutive for $d_1,...,d_{n-r}$ and thus also for $d_1,...,d_n$ by adopting the ordering $(d_{n-r+1},...,d_n,d_1,...,d_{n-r})$. It allows to define a differential module $M_r$ and an epimorphism $M_r\rightarrow M\rightarrow 0$ as $M$ {\it is defined by more equations}. Now, as $cd_D(M)=r$, we have $t_{r-1}(M)=M$ and each torsion element of $t(M)=M$ surely satisfies {\it at least} $r$ PD equations involving successively $d_n,...,d_{n-r+1}$. As for the other equations, they should {\it only} include $d_1,...,d_{n-r}$ and this is just the way to construct $t(M_r)$ by considering the exact sequence:\\
\[  0 \longrightarrow t(M_r) \longrightarrow M_r \longrightarrow k({\chi}_1,...,{\chi}_{n-r})\otimes M_r  \]
As it is well known that localization preserves exactness [3,22], this exact sequence projects onto the desired one.\\
\hspace*{12cm}   Q.E.D. \\

\noindent
{\bf REMARK 3.9}: Using modules instead of ideals, the above theorem allows to generalize for arbitrary $m$ the condition obtained by Macaulay for $m=1$ ([13], \S 41, p 43 and \S 43, p 45) that we translate into modern language as another proof (See [3], IV, \S 1,  exercise 9). In fact, if $\mathfrak{a}\subset A=k[{\chi}_1,...,{\chi}_n]$ is such that $cd(A/ \mathfrak{a})=r$, then $\mathfrak{a}$ is unmixed ($A/\mathfrak{a}$ is r-pure) if and only if $S({\chi}_1,...,{\chi}_{n-r})P\in \mathfrak{a}\Rightarrow P\in \mathfrak{a}$. For if $\mathfrak{a}={\mathfrak{q}}_1\cap ... \cap {\mathfrak{q}}_s$ with, say $cd(A/ {\mathfrak{q}}_1)>r$, then $\exists S\in {\mathfrak{q}}_1$ and $P\notin {\mathfrak{q}}_1, P\in {\mathfrak{q}}_2\cap ... \cap {\mathfrak{q}}_s$ so that $SP\in \mathfrak{a}$ does not require $P\in \mathfrak{a}$. Conversely, if no such ${\mathfrak{q}}_1$ exists and $SP\in \mathfrak{a}$, then $P\in \mathfrak{a}:(S)=\mathfrak{a}$. Indeed and more generally, if $\mathfrak{b}\in A$ is an ideal and $\mathfrak{a}\subset \mathfrak{a}:\mathfrak{b}\neq \mathfrak{a}$, then $\mathfrak{b}\subset {\mathfrak{p}}_i=rad({\mathfrak{q}}_i)$ for some $i$. To prove this, if $\mathfrak{a}:\mathfrak{b}=\mathfrak{c}$ and $\mathfrak{b}$ is not contained in ${\mathfrak{p}}_i$, then one can find $b\in \mathfrak{b}, b\notin {\mathfrak{p}}_i $ and $c\in \mathfrak{c}, c\notin \mathfrak{a}$ with $bc\in \mathfrak{a}\subset {\mathfrak{q}}_i$ and thus $c\in {\mathfrak{q}}_i,\forall i$ that is $c\in \mathfrak{a}$ and a contradiction.\\

The two following corollaries generalize Proposition 2.9:\\

\noindent
{\bf COROLLARY 3.10}: ([13], \S 41, p 43) $M$ is r-pure $\Leftrightarrow  0\rightarrow M \rightarrow k({\chi}_1,...,{\chi}_{n-r})\otimes M$ is exact.\\

\noindent
{\bf COROLLARY 3.11}: $M$ is r-pure $\Leftrightarrow \exists 0\rightarrow M \rightarrow L$ with {\it projective dimension} $pd_D(L)=r$ if $r\geq 1$.\\
{\it Proof}: (Compare to [2], p 494 and [18], p 553) As $M$ is r-pure, then $t_r(M)=0$ that is $M_r$ is torsion free and we may use $k({\chi}_1,...,{\chi}_{n-r})\otimes M_r$ in order to parametrize $M_r$ exactly as we did for embedding a torsion-free module into a free module. According to Proposition 2.16, the parametrization now depends on ${\alpha}^{n-r}_1$ arbitrary unknowns $z$, that is we may embed $M_r$ into ${\alpha}^{n-r}_1$ copies of $k[d_1,...,d_{n-r}]$ when coming back to the differential framework. After substitution into the original equations, the equations of class 1 up to class $(n-r)$ disappears for the $z$ as they are {\it automatically} satisfied by the parametrization. The number of nonmuliplicative derivatives is $\leq r-1$ ({\it care}) for each of the remaining equations of class $(n-r+1)$ up to class $n$. But such a number is just the way to know about the projective/free dimension by constructing a resolution of $M$ or a Janet sequence for the system ([17], p 146).\\
\hspace*{12cm}   Q.E.D. \\

\noindent
{\bf REMARK 3.12}: In actual practice it is important to notice that the partial localization {\it kills} the PD equations of class 1 up to class $(n-r-1)$ ({\it care again}) because of the compatibility conditions provided by the involutive assumption. Moreover we now understand why Macaulay ([13], \S 79, p 89) was always dealing with unmixed ideals $\mathfrak{a}$ or pure modules $A/\mathfrak{a}$. It is known ([2] and [18], Proposition 3.173, p 549) that an $A$-module $M$ is r-pure if and only if $cd(A/\mathfrak{p})=cd(M)=r$ for any $\mathfrak{p}\in spec(A)$ appearing in the prime decomposition of $rad(ann_A(M))$ {\it and} no embedded primary components occur in a primary decomposition of $ann_A(M)$ or equivalently $ass(M)=\{\mathfrak{p}\in spec (A)\mid \exists x\in M, ann(x)=\mathfrak{p}\}$ is equidimensional. Example 2.3 shows that the second condition is needed. Accordingly, {\it any prime or primary module is pure}.\\

\noindent
{\bf EXAMPLE 3.13}: $k=\mathbb{Q}, m=1,n=3,q=2,r=1$. The module $M$ defined by the involutive second order system ${\Phi}^3\equiv y_{33}=0, {\Phi}^2\equiv y_{23}=0, {\Phi}^1\equiv y_{13}=0$ is not pure. Among the three compatibility conditions we have $d_2{\Phi}^1-d_1{\Phi}^2=0$. As only the class 3 is full, the localization is done by tensoring with $k({\chi}_1,{\chi}_2)$ and we get ${\Phi}^1=\frac{{\chi}_1}{{\chi}_2}{\Phi}^2$. Also, as $y_1$ and $y_2$ are killed by $d_3$ they are in $t_0(M)=M$ but not in $t_1(M)$. Hence there is a {\it gap} because $t_1(M)=t_2(M)$ and $y_3$ generates $t_2(M)$ as it is killed by $(d_1,d_2,d_3)$ and $t_3(M)=0$ by definition. In order to get a first order presentation (though {\it not} a Spencer form) we may introduce $z^1=y,z^2=y_1,z^3=y_2, z^4=y_3$ and $M_1$ is defined by 3 equations of class 2 and 2 equations of class 1. With respect to $(d_1,d_2)$, $M_1$ is {\it not} torsion-free but {\it not} a torsion module as $t(M_1)$ is generated by $z^4$ as already noticed. Finally we have $ass(M)=\{({\chi}_3),({\chi}_1,{\chi}_2,{\chi}_3)\} $as another way to check that $M$ is not pure.\\

\noindent
{\bf EXAMPLE 3.14}: $k=\mathbb{Q}, m=3, n=4, q=1, r=1$. The module $M$ defined by the first order homogeneous involutive system $y^1_4=0,y^2_4=0,y^3_4=0,y^3_3+y^2_2+y^1_1=0$ is 1-pure. We notice that $M_1$ is defined by the only divergence-free condition and is thus torsion-free. Indeed, tensoring by $k({\chi}_1,{\chi}_2,{\chi}_3)$ in order to localize, we get the parametrization $y^3=-\frac{{\chi}_2}{{\chi}_3}y^2-\frac{{\chi}_1}{{\chi}_3}y^1=-{\chi}_2z^2-{\chi}_1z^1, y^2={\chi}_3z^2,y^1={\chi}_3z^1$ and we have a strict embedding $M\subset L$ with $L$ generated by $(z^1,z^2)$ satisfying only $z^1_4=0,z^2_4=0$. Accordingly, $L$ admits a resolution $0\rightarrow D^2 \rightarrow D^2 \rightarrow L \rightarrow 0$ with morphism $(P_1,P_2)\rightarrow (P_1d_4,P_2d_4)$ and $pd(L)=1$.\\

\noindent
{\bf EXAMPLE 3.15}: With $m=3,n=1, q=2$, the module defined by ${\Phi}^1\equiv y_{33}=0, {\Phi}^2\equiv y_{13}-y_2=0$ (Example 2.13) and the module defined by ${\Phi}^1\equiv y_{33}-y_3=0, {\Phi}^2\equiv y_{13}-y_2=0$ are 2-pure and have already a projective dimension equal to 2. Indeed, using computer algebra as in ([7],[21]), then $({\Phi}^1, {\Phi}^2)$ does satisfy a single second order CC in both cases.\\

\noindent
4) {\bf INVERSE SYSTEMS}\\

Let $K$ be a differential field with subfield of constants $k=cst(K)$. The ring $D=K[d]$ is filtred by the order $q$ of operators and we have $K=D_0\subset D_1\subset ... \subset D_{\infty}=D$. Accordingly, as explained at the end of section 2, the module $M$ is filtred by the order $q$ of the linear combinations $y_q=D_qy$ allowing to describe elements of $M$ by residue and we have the {\it inductive limit} $M_0\subseteq M_1\subseteq ... \subseteq M_q\subseteq ... \subseteq M_{\infty}=M$ with $d_iM_q\subseteq M_{q+1}$ and $M=DM_q$ for $q\gg 0$. For example, according to the last example where $k={\mathbb{Q}}$, the system $y_{33}=0, y_{13}-y_2=0$ defines a module $M$ having the finite free resolution $0\rightarrow D \rightarrow D^2 \rightarrow D \rightarrow M \rightarrow 0$. In order to determine $M_2$, one has to take the residue of $D_2$ with respect to the vector space 
$k y_{33} + k y_{23} + k y_{22} + k (y_{13}-y_2)$ which is the intersection of $D_2$ with the image of the presentation morphism $D^2 \rightarrow D: (P,Q) \rightarrow P y_{33}+Q(y_{13}-y_2)$ which is not strict, a result showing that it is important to start with an involutive operator or at least with a strict presentation.\\

\noindent
{\bf DEFINITION 4.1}: We define the {\it system} $R=hom_K(M,K)=M^*$ and set $R_q=hom_K(M_q,K)=M_q^*$ as the {\it system of order q} in order to have now the {\it projective limit} $R=R_{\infty}\rightarrow ... \rightarrow R_q \rightarrow ... \rightarrow R_1 \rightarrow R_0$. Taking into account the differential geometric framework of Spencer ([17],[18],[24]), if a system of PD equations of order $q$ is given as before, then $f_q\in R_q:y^k_{\mu}\rightarrow f^k_{\mu}\in K$ with $a^{\tau\mu}_kf^k_{\mu}=0$ defines a {\it section at order} $q$ and we set $f_{\infty}=f\in R$ for a {\it section}. It is only when the field of constants $k$ is used that we can speak about a formal power series solution (see Example 4.5 and all the examples of section 5 for explicit finite or infinite situations).\\

\noindent
{\bf REMARK 4.2}: In the case of an involutive system of order $q$ in {\it solved form}, the matrix $(a^{\tau\mu}_k)$ and the corresponding prolongations for increasing orders allow to express certain jets, 
called {\it principal}, from the other jets, called {\it parametric} (called "{\it complete set of remainders}" by Macaulay), and exacty describe the upper part of the picture drawn by Macaulay in ([13], \S 59, p 67 and \S 68, p 79). Similarly, for any $q\geq 0$ the following commutative and exact diagram:\\

\[ \begin{array}{rcccccccl}
   &      &                     & 0      &                     & 0             &                     &          &  \\
   &      &                     &\downarrow &         &\downarrow &               &          &  \\
   &  0  & \rightarrow & R^q & \rightarrow & R^{q-1} & \rightarrow &  g_q & \rightarrow 0\\
   &      &                      &\downarrow &        &\downarrow &               &          &  \\
   &  0  & \rightarrow & R     &  =                 &  R           & \rightarrow &  0      &  \\
   &       &                     &\downarrow &        &\downarrow &               &          &   \\
0\rightarrow& g_q & \rightarrow & R_q &\stackrel{{\pi}^q_{q-1}}{\rightarrow} & R_{q-1} &\rightarrow &0& \\
   &        &                     &\downarrow &     & \downarrow &             &             &   \\
   &        &                     &    0                &     &    0               &                &           &    
   \end{array}  \]

\noindent   
allows to define the {\it symbol} $g_q$ of $R_q$ and the upper row is again {\it exactly} describing the lower part of the same picture through the use of {\it truncated} formal power series or sections. The symbols are a modern way to describe the {\it compartments} of Macaulay. The use of a basis $(1,0,....), (0,1,....) $ and so on for the parametric jets, after ordering them, brings a block triangular matrix as explained by Macaulay ([13], \S 59, p 67).\\

\noindent
{\bf DEFINITION 4.3}: A {\it modular equation} $E\equiv a^{\mu}_kf^k_{\mu}=0$ of order $q$ with $0\leq k\leq m, 0\leq \mid\mu\mid\leq q$ is just a way to write down a section $f_q\in R_q$ by using an implicit summation with {\it formal coefficients}. Of course, as in Example 4.5 below, infinite summations may also be considered. The procedure is absolutely similar to the case $m=1,K=k$ where one uses the purely formal power series notation $\sum f_{\mu}\frac{x^{\mu}}{\mu !}$ for writing down a section, even though the variable $x$ has absolutely no meaning in the module framework. Finally, as noticed by Macaulay, if one considers the set of modular equations at order $q$ as a homogeneous linear system for the unknowns $a^{\mu}_k$ at order $q$, then of course the given coefficients $a^{\tau\mu}_k$ form a basis of solutions linearly independant over $K$ and indexed by $\tau$. This is the reason for which we have chosen a similar notation.\\

The following proposition generalizes the results of Macaulay to arbitrary systems with variable coefficients because $K$ is a $D$-module with the standard action $(D,K)\rightarrow K: (d_i,a)\rightarrow {\partial}_ia$. However, it is not evident, at first sight, to endow $M^*$ with a structure of left $D$-module in general, unless $D$ is a commutative ring, that is $K=cst(K)=k$ ([2], Theorem 1.3.1, 21, [18], Theorem 3.89, p 487).\\

\noindent
{\bf PROPOSITION 4.4}: When $M$ is a left $D$-module, then $R$ is a left $D$-module too.\\
\noindent
{\it Proof}: It is clear that $D$, as an algebra, is generated by $K=D_0$ and $T=D_1/D_0$ with $D_1=K\oplus T$. Let us define:\\
\[   (af)(m)=af(m)=f(am) \hspace{1cm} \forall a\in K, \forall m\in M\]
\[   (\xi f)(m)=\xi f(m)-f(\xi m)  \hspace{1cm}  \forall \xi =a^id_i\in T, \forall m\in M  \]
It is easy to check that $d_ia=ad_i+{\partial}_ia$ in the operator sense and that $\xi\eta -\eta\xi =[\xi,\eta]$ is the standard bracket of vector fields. We finally get $(d_if)^k_{\mu}=(d_if)(y^k_{\mu})={\partial}_if^k_{\mu}-f^k_{\mu +1_i}$ that is {\it exactly} the {\it Spencer operator} ([17],[18],[24]). We have $d_id_j=d_jd_i=d_{ij}, \forall i,j=1,...,n$ because $(d_id_jf)^k_{\mu}={\partial}_{ij}f^k_{\mu}-{\partial}_if^k_{\mu +1_j}-{\partial}_jf^k_{\mu +1_i} +f^k_{\mu +{1_i}+{1_j}}$ and $d_iR_{q+1}\subseteq R_q$, a result leading to $d_iR\subset R$  and a well defined operator $R\rightarrow T^*{\otimes}_K R:f \rightarrow dx^i {\otimes} d_if$. This is the dual framework of the Spencer resolution $D{\otimes}_K T{\otimes}_K M \rightarrow D {\otimes}_K M \rightarrow M \rightarrow 0$ with $P\otimes \xi \otimes m \rightarrow P\xi\otimes m - P\otimes \xi m$ and $P\otimes m \rightarrow Pm$ ([2],[11], p 1,19, [18], p 499).\\
Alternatively and in a coherent way with differential geometry, if we have a linear system ${\Phi}^{\tau}=0$ defining $R_q$ and its first prolongation ${\Phi}^{\tau}=0,d_i{\Phi}^{\tau}=0 $ defining $R_{q+1}$, as already exhibited in Section 2, a section $f_{q+1}\in R_{q+1}$ over $f_q\in R_q$ satisfies both $a^{\tau\mu}_kf^k_{\mu}=0$ and $a^{\tau\mu}_kf^k_{\mu +1_i}+{\partial}_ia^{\tau\mu}_kf^k_{\mu}=0$ as equalities in $K$ with $0\leq \mid\mu\mid \leq q$. Applying ${\partial}_i$ to the first and substracting the second, we get $a^{\tau\mu}_k({\partial}_if^k_{\mu}-f^k_{\mu +1_i})=0$. Accordingly, we obtain:\\
\[  f_{q+1}\in R_{q+1}\stackrel{d_i}{\longrightarrow}d_if_{q+1}\in R_q\Longleftrightarrow 
E\equiv a^{\nu}_kf^k_{\nu}=0  \stackrel{d_i}{\longrightarrow} d_iE\equiv a^{\mu}_k({\partial}_if^k_{\mu}-f^k_{\mu +1_i})=0, \forall f\in R   \]
but $d_iE$ is of order $q$ with $0\leq\mid\mu\mid\leq q$ whenever $E$ is of order $q+1$ with $0\leq\mid\nu\mid\leq q+1$. When $K=k$, the partial derivative disappears and we recognize, {\it exactly but up to sign},  the operator of Macaulay ([13], \S 60, p 69). For this reason and unless mentioned explicitly, in this specific 
situation we shall change the sign of the Spencer operator in order to agree with Macaulay.\\
\hspace*{12cm}   Q.E.D.  \\

As $D$ is a left $D$-module, it follows from the above proposition that $D^*=hom_K(D,K)$ is a left $D$-module. Moreover, it is known ([4], Proposision 11, p 18,[15],[22], p 37) that $D^*$ is an injective $D$-module as there is a canonical isomorphism $M^*=hom_K(M,K)\simeq hom_D(M,D^*)$ where both sides are well defined. It is also worth pointing out the importance of the two preceding propositions for computer algebra as they allow to deal with sections and {\it not} with solutions, contrary to the tradition. We emphasize once more that this new point of view is the main tool brought by Spencer and leading to the {\it Spencer sequences} ([17],[23]). Indeed, starting with a presentation $Dz \rightarrow Dy: {z}^{\tau}_{\nu} \rightarrow d_{\nu}{\Phi}^{\tau}$ as previously defined and introducing an arbitrary section $y^k_{\mu} \rightarrow {\xi}^k_{\mu}\in K$ leading by composition to a section $z^{\tau}_{\nu} \rightarrow {\eta}^{\tau}_{\nu}\in K$, we obtain from the proof of the last proposition ${\partial}_i{\eta}^{\tau} - {\eta}^{\tau}_i=a^{\tau\mu}_k ({\partial}_i{\xi}^k_{\mu}-{\xi}^k_{\mu +1_i})$ a result showing that the Spencer operator commutes with the dual of the presentation ([17], p 147, [23]). \\

\noindent
{\bf EXAMPLE 4.5}: $k=\mathbb{Q}, m=1, n=1, q=2$. For the system $y_{11}-y=0$, if we set:\\
$   f'=(1,0,1,0,...)\rightarrow E' \equiv a^0+a^{11}+ ... =0$ and $f''=(0,1,0,1,...)\rightarrow E" \equiv a^1+a^{111}+ ... =0$\\
and use $\{f',f''\}$ as a basis of $R$ over $k$, we have $d_1f''=f'$ or equivalently $d_1E''=E'$.\\

\noindent
{\bf REMARK 4.6}: If $cd(M)=r$ and $K=k$, then ${\alpha}^{n-r+1}_1=0,..., {\alpha}^n_1=0$ and a partial localization  brings the system to a {\it finite type} (zero symbol at high order) system in $(d_{n-r+1},...,d_n)$ over the field $k({\chi}_1,...,{\chi}_{n-r})$. Accordingly, there is a finite number of linearly independent corresponding sections of the localized system and thus an equal finite number of (dual) modular equations as in the previous example ([13], \S 79, p 88). In the situation $K=k$, we have also $ann_D(M)=ann_D(R)$. Indeed, as a representative of any element of $M$ can be written as a {\it finite} linear combination of parametric jets with coefficients in $k$, we have $M\subseteq M^{**}$ and thus $ann(M)\subseteq ann(M^*)\subseteq ann(M^{**})\subseteq ann(M) \Rightarrow ann(M)=ann(M^*)$. This result generalizes the one of Macaulay ([13], \S 61, p 70) obtained when $m=1$. Indeed, $a^{\tau\mu}f_{\mu}=0\Rightarrow a^{\tau\mu}f_{\mu +\nu}=0$ by prolongation and thus $E\equiv a^{\nu}f_{\nu}=0\Rightarrow a^{\tau\mu}d_{\mu}E\equiv a^{\tau\mu}a^{\nu}f_{\mu +\nu}=0$. We may also set $d_iE\equiv a^{\nu-1_i}f_{\nu}$ with $a^{\nu-1_i}=0 $ if ${\nu}_i=0$.\\

Following ([18], p 113), any primary decomposition, say with two components for simplicity, gives rise to a monomorphism $0\rightarrow M\rightarrow Q' \oplus Q''$ where $Q',Q''$ are primary modules, both with two epimorphisms $M\rightarrow Q' \rightarrow 0, M\rightarrow Q'' \rightarrow 0$, respectively induced by the localization morphisms $M\rightarrow M_{{\mathfrak{p}}{'}}, M\rightarrow M_{{\mathfrak{p}}{"}}$ when $M$ is pure (unmixed annihilator) with $ass(Q')=\{{\mathfrak{p}}'\}, ass(Q")=\{{\mathfrak{p}}"\}$ and $ass(M)=\{{\mathfrak{p}}{'},{\mathfrak{p}}{"}\}$. Setting $R'=hom_K(Q',K), R''=hom_K(Q'',K)$ and using the fact that $D^*$ is injective, we get an epimorphism $R' \oplus R'' \rightarrow R\rightarrow 0$ both with two monomorphisms $0\rightarrow R' \rightarrow R, 0\rightarrow R'' \rightarrow R$ proving that $R' , R'' , R' + R'' , R' \cap R'' $ are subsystems of $R$. The following proposition, {\it not evident at first sight}, explains the aim of Macaulay ([13], end of \S 79, p 89) and allows to {\it use various subsystems for studying} $R$ {\it instead of decomposing} $M$.\\

\noindent
{\bf PROPOSITION 4.7}:   $ R = R' + R'' $.\\
\noindent
{\it Proof}: We have the well known short exact sequence  $0 \rightarrow R' \cap R'' \rightarrow R' \oplus R'' \rightarrow R' + R'' \rightarrow 0 $ where the last morphism is $(f',f'')\rightarrow (f'-f'')$. Composing the epimorphism with the monomorphism $0 \rightarrow R' + R'' \rightarrow R $ and using the fact that the composite morphism $ R' \oplus R'' \rightarrow R$ is an epimorphism, it follows that the previous monomorphism is also an epimorphism and thus an isomorphism, though in general $R' \cap R'' \neq 0$, unless we have ${\mathfrak{p}}' + {\mathfrak{p}}'' = A$, a situation always met with $max(A)\subset spec(A)$ in the case of modules over a ring $A$ which is also a finitely generated algebra over a field $k$.\\
\hspace*{12cm}     Q.E.D.  \\

We finally recall in a self-contained way a few results on the so-called {\it socle} and {\it top} of a module $M$ over a commutative noetherian integral domain $A$ with unit $1$ ([1],Ê[6]). First of all, we quote the following theorem on associated primes where {\it both} isolated {\it and} embedded components are needed ([3], IV, \S 1, exercise 11).\\

\noindent
{\bf THEOREM 4.8}: If $M$ is a finitely generated $A$-module, the sequence $0 \rightarrow M \rightarrow {\oplus}_{\mathfrak{p}\in ass(M)}M_{\mathfrak{p}}$ is exact.\\
\noindent
{\it Proof}: If the sequence is not exact, let $N$ be the kernel of the morphism on the right. If $ass(M)=\{{\mathfrak{p}}_1, ... ,{\mathfrak{p}}_s\}$, let us consider the defining exact sequences $0\rightarrow N_i\rightarrow M \rightarrow M_{{\mathfrak{p}}_i}, \forall i=1,...,s$. By definition, we have $N=\cap N_i$ and it is well known that $N\neq 0 \Leftrightarrow ass(N)\neq \emptyset$. In that case, let $\mathfrak{p}\in ass(N)\subset ass(M)$, that is to say $\mathfrak{p}={\mathfrak{p}}_i$ for a certain $1\leq i\leq s$. Again by definition, one can find $x\in N\subset N_i$ such that $\mathfrak{p}={\mathfrak{p}}_i=ann(x)$. But $x\in N_i \Leftrightarrow \exists s_i\in S_i=A-{\mathfrak{p}}_i, s_ix=0$ because of localization and thus a contradiction. One could also say that $ann(x)\subset {\mathfrak{p}}_i$ for some $i$ whenever $0\neq x\in N\subset M$ as it is well known that $\cup {\mathfrak{p}}_i$ is the set of zero divisors for $M$. But $x\in N_i$ and we conclude as above.\\
\hspace*{12cm}          Q.E.D.\\

\noindent
{\bf REMARK 4.9}: When $t(M)\neq M$, then $(0)\in ass(M)$ and the image of the canonical morphism $M \rightarrow M_{(0)}$ is just $M/t(M)$ as in Remark 2.10. However, we recall that one can embed a module into a direct sum of primary modules by using the images of the morphisms $M\rightarrow M_{\mathfrak{p}}$ in the preceding theorem on the condition that $ass(M)$ only contains minimal primes ([3],[18]). Such a situation happens in the case of pure modules or in the case of quotients of unmixed ideals considered by Macaulay (se Examples 2.2 to 2.4 and 3.15).\\

\noindent
{\bf EXAMPLE 4.10}: Let $A=\mathbb{Q}[x,y,z], M=A/{\mathfrak{a}}$ with $\mathfrak{a}=(x^2,xy,xz,yz)={\mathfrak{p}}_1\cap{\mathfrak{p}}_2\cap {\mathfrak{m}}^2$ where ${\mathfrak{p}}_1=(x,y), {\mathfrak{p}}_2=(x,z)$ are the two minimal primes (isolated components of the characteristic variety) and $ \mathfrak{m}=(x,y,z)\in max(A)$ (embedded component). Then $ass(M)=\{ {\mathfrak{p}}_1, {\mathfrak{p}}_2, \mathfrak{m}\}$ where ${\mathfrak{p}}_1$ kills $\bar{z}$, ${\mathfrak{p}}_2$ kills $\bar{y}$ and $\mathfrak{m}$ kills $\bar{x}$. It follows that $\bar{x}$ belongs to the kernel of $M\rightarrow M_{{\mathfrak{p}}_1}\oplus M_{{\mathfrak{p}}_2}$ and cannot be killed by any $s\in A-\mathfrak{m}$.\\

Keeping in mind the bricks needed in order to construct a house, {\it a basic idea in module theory is to look for the greatest semi-simple submodule of a given module}. For this, if $\mathfrak{m}\in max(A)\cap ass(M)$, then one can find a finite number of elements $x,y,...\in M$ killed by $\mathfrak{m}$. Accordingly, the map $x:A\rightarrow M:a\rightarrow ax$ has kernel $\mathfrak{m}$ and $A/\mathfrak{m}\simeq Ax\subseteq M$ is a simple module, like $Ay$ which may eventually be different and so on. The direct sum $Ax\oplus Ay\oplus ...$ is called the {\it socle} of $M$ at $\mathfrak{m}$ and denoted by ${soc}_{\mathfrak{m}}(M)$. These simple components are called {\it isotypical} as they are all isomorphic to $A/\mathfrak{m}$.\\

\noindent
{\bf DEFINITION 4.11}: The {\it socle} of $M$ is $soc(M)=\oplus {soc}_{\mathfrak{m}}(M)$ for $\mathfrak{m}\in  max(A)\cap ass(M)$. It is the largest semi-simple submodule of $M$. \\

We notice that the double condition on the direct sum is essential as we need not only a submodule ($\mathfrak{m}\in ass(M)$) but {\it also} a simple module ($\mathfrak{m}\in max(A)\subseteq spec(A)$). Also, if $S{'}, S{"}$ are two simple submodules of $M$, then $M/S{'}\oplus S{"}$ is the {\it fiber sum} of $M/S{'}$ and $M/S{"}$ over $MÊ$ ([12], p 88) but the resulting construction is not natural and will provide a motivation for duality in order to use Proposition 4.7. Finally, $M$ is semi-simple if $M=soc(M)$ and $soc(M)=0$ if $M$ has no simple submodule, like the $\mathbb{Z}$-module $\mathbb{Z}$. In the previous example $S=A\bar{x}\simeq A/\mathfrak{m}$ is the only simple submodule of $M$.\\

\noindent
{\bf EXAMPLE 4.12}: If $A=\mathbb{Q}[x,y]$ and $M=A/\mathfrak{a}$ with $\mathfrak{a}=(x^3,y^2,xy)$, then {\it both} $\bar{y}$ and ${\bar{x}}^2$ are killed by $\mathfrak{m}=(x,y)$. It follows that $soc(M)=soc_{\mathfrak{m}}(M)=A\bar{y}\oplus A{\bar{x}}^2$ has two isotypical components isomorphic to $A/\mathfrak{m}$.\\

\noindent
{\bf LEMMA 4.13}: Any morphism $f:M\rightarrow N$ induces a morphism $f:soc(M)\rightarrow soc(N)$. In particular, if $M{'}$ is a submodule of $M$, then $soc(M{'})=M{'}\cap soc(M)$.\\
{\it Proof}: The lemma follows at once from the Schur lemma saying that, when $f\neq 0$, then $M$ simple $\Rightarrow f$ injective, $N$ simple $\Rightarrow f$ surjective.\\
\hspace*{12cm}      Q.E.D.   \\

\noindent
{\bf DEFINITION 4.14}: The {\it radical} of a module $M$ is the submodule $rad(M)$ which is the intersection of all the maximum {\it proper} submodules of $M$. If $rad(M)=0$, for example if $M$ is simple, we say that $M$ has no radical. If $M$ has no proper maximum submodule, then $rad(M)=M$.\\

\noindent
{\bf LEMMA 4.15}: $rad(M)$ is the intersection of all the kernels of the nonzero morphisms $M\rightarrow S$ where $S$ is a simple module.\\
\noindent
{\it Proof}: From the Schur lemma, the above morphism is an epimorphism and we may introduce the defining short exact sequences $0 \rightarrow N \rightarrow M \rightarrow S \rightarrow 0$. Let us consider the exact sequence $0\rightarrow \cap N \rightarrow M \rightarrow \oplus S $. The image of the morphism on the right is a submodule of a semi-simple module and thus a semi-simple module too, which is even a direct summand. Accordingly, restricting the choice of the simple modules in order to have an irredundant intersection, the morphism on the right thus becomes an epimorphism leading to the next definition.\\
\hspace*{12cm}    Q.E.D. \\

\noindent
{\bf DEFINITION 4.16}: The {\it top} of the module $M$ is the semi-simple module defined by the short exact sequence $0\rightarrow rad(M)\rightarrow M \rightarrow top(M)\rightarrow 0$. It can also be defined as the largest quotient of $M$ that is a direct sum of simple modules.\\

We have the following three useful lemmas [1]:\\

\noindent
{\bf LEMMA 4.17}: Any morphism $f:M\rightarrow N$ induces a morphism $f=rad(M)\rightarrow rad(N)$.\\
\noindent
{\it Proof}: Let $S$ be a simple module. For any morphism $g: N \rightarrow S$, the composition $g\circ f:M \rightarrow S$ vanishes on $rad(M)$ and thus $g$ vanishes on $f(rad(M))$, that is $f(rad(M))\subseteq rad(N)$.\\
\hspace*{12cm}       Q.E.D.  \\

\noindent
{\bf LEMMA 4.18}: If $M\neq 0$ is finitely generated, then $rad(M)\neq M$.\\
\noindent
{\it Proof}: From noetherian arguments, $M$ always contains a maximum proper submodule.\\
\hspace*{12cm}   Q.E.D.   \\

\noindent
{\bf LEMMA 4.19}: (Nakayama) Let $M$ be a finitely generated module and $N$ a submodule of $rad(M)$. If $L\subseteq M$ is such that $L+N=M$, then $L=M$.\\
\noindent
{\it Proof}: Let us suppose that $L\neq M$. Then, from noetherian arguments again, $L$ is contained in a maximum proper submodule $L{'}$. It follows that $N+L\subseteq rad(M)+L{'}\subseteq L{'}$ and a contradiction.\\
\hspace*{12cm}   Q.E.D.  \\

We are now ready to provide the achievement of this paper while explaining ([13], \S 77,79,82).\\

\noindent
5) {\bf MACAULAY ' S SECRET}\\

{\it The crucial idea of Macaulay has been to use} $top(R)$ {\it instead of} $soc(M)$ {\it by means of duality theory}, in order to use Nakayama's lemma for finding generating sections (formal solutions) of $R$. We proceed in a few successive steps for working with differential modules in an effective way and treating the following specific examples.\\

5.1) The first basic procedure is to check that $M$ is r-pure. For this we {\it must} determine $r$ by exhibiting an involutive system. As already noticed, a linear change of derivations may be needed in order to check involution.The partial localization will then be used in order to check the purity and to deal only with maximal ideals because a prime ideal is maximum if and only if its residue integral domain is zero dimensional.\\

\noindent
{\bf EXAMPLE 5.1.1}: If $\mathfrak{a}=({\chi}_1,{\chi}_2{\chi}_3)=({\chi}_1,{\chi}_2)\cap ({\chi}_1,{\chi}_3)$ the corresponding system $y_1=0, y_{23}=0$ is not involutive and the change $ {\chi}_1\rightarrow {\chi}_3, {\chi}_2\rightarrow {\chi}_2, {\chi}_3\rightarrow {\chi}_2-{\chi}_1$ provides the involutive system in solved form $y_{33}=0, y_{23}=0, y_{22}-y_{12}=0, y_{13}=0, y_3=0$. We have $M\subset \mathbb{Q} ({\chi}_1)\otimes M$ and the localized module has the two associated maximum ideals ${\mathfrak{m}}_1=(d_3,d_2)$ and ${\mathfrak{m}}_2=(d_3,d_2-{\chi}_1)$ with 
${\mathfrak{m}}_1 + {\mathfrak{m}}_2=\mathbb{Q}({\chi}_1)[d_2,d_3]$ as $d_2 - (d_2-{\chi}_1)={\chi}_1$. \\

\noindent
{\bf EXAMPLE 5.1.2}: $\mathfrak{a}=(({\chi}_1)^3,({\chi}_2)^2,{\chi}_1{\chi}_2)$ is primary because $rad(\mathfrak{a})=({\chi}_1,{\chi}_2)=\mathfrak{m}$. With $D=\mathbb{Q}[d_1,d_2]$, the homogeneous second order system $R_3=\{ y_{222}=0, y_{122}=0, y_{112}=0, y_{111}=0, y_{22}=0, y_{12}=0\}$ is trivially involutive because its symbol is zero. The corresponding module is primary and $2$-pure. No localization is needed and $dim_{\mathbb{Q}}(R)=dim_{\mathbb{Q}}(R_3)=4$.\\

\noindent
{\bf EXAMPLE 5.1.3}: $\mathfrak{a}=(({\chi}_3)^2, {\chi}_2{\chi}_3-({\chi}_1)^2, ({\chi}_2)^2)$ is primary because $rad(\mathfrak{a})=({\chi}_1,{\chi}_2,{\chi}_3)=\mathfrak{m}$. With now $D=\mathbb{Q}[d_1,d_2,d_3]$, the homogeneous system $R_2=\{ y_{33}=0, y_{23}-y_{11}=0, y_{22}=0\}$ is {\it not} involutive (see [18, p 321 for another similar example) but its prolongation $R_4$ is trivially involutive 
with symbol $g_4=0$ and $dim_{\mathbb{Q}}(g_3)=1$. No localization is needed and $dim_{\mathbb{Q}}(R)=dim_{\mathbb{Q}}(R_2)=8$.\\

\noindent
5.2)  The idea is now to adapt to modules an argument already used in Remark 3.9 for ideals. If $\mathfrak{a}=\cap {\mathfrak{q}}_i$ is a primary decomposition with ${\mathfrak{p}}_i=rad({\mathfrak{q}}_i)$, then $\mathfrak{a}:\mathfrak{b}\neq \mathfrak{a} \Rightarrow {\mathfrak{b}}\subset {\mathfrak{p}}_i$ for a certain $i$. In particular, $\mathfrak{a}:\mathfrak{m}\neq \mathfrak{a}$ for $\mathfrak{m}\in max(A)\Rightarrow \mathfrak{m}={\mathfrak{p}}_i$ for a certain $i$.\\

Let us set $A^m=F$ as a free module and consider the short exact sequence $0\rightarrow I \rightarrow F \rightarrow M \rightarrow 0$ where $I$ is the so-called {\it module of equations} and let $\mathfrak{a}\in A$ be an ideal. We want to prove the following lemma:\\

\noindent
{\bf LEMMA 5.2.1}: $I\subseteq I:\mathfrak{a}=J\neq I \Rightarrow \mathfrak{a}\subseteq {\mathfrak{p}}_i$ for a certain $i$.\\
\noindent
{\it Proof}: Let us consider a primary decomposition $I=\cap I_i$ in $F$ and pass to the quotient by introducing short exact sequences $0\rightarrow I_i\rightarrow F\rightarrow Q_i\rightarrow 0$ in order to have epimorphisms $M\rightarrow Q_i\rightarrow 0$ and a monomorphism $0\rightarrow M\rightarrow \oplus Q_i$. If $I\neq J$, let $x\in J, x\notin I$ with $ax\in I\Rightarrow a\bar{x}=0, \bar{x}\neq 0, \forall a\in \mathfrak{a}$. It follows that $a$ is a zerodivisor and thus $a\in \cup {\mathfrak{p}}_i$. Also, $x\notin I_i$ for a certain $i$ otherwise $x\in \cap I_i=I$. But $I_i$ is (co)primary with $ax\in I\subset I_i,x\notin I_i $ (or $a\bar{x}=0, \bar{x}\neq 0$ in $Q_i$) $\Rightarrow a\in {\mathfrak{p}}_i\Rightarrow \mathfrak{a}\subset {\mathfrak{p}}_i$. In particular, if $I:\mathfrak{m}\neq I$ for $\mathfrak{m}\in max(A)$, then $\mathfrak{m}={\mathfrak{p}}_i$ for a certain $i$ as before and $\mathfrak{m}\in ass(M)$.\\
\hspace*{12cm}    Q.E.D.   \\

It is at this precise point that we have to use specific properties of the ring $D$ that will now be used in place of $A$. From now on and unless specified, we assume that the partial localization has been realized. Using therefore $k({\chi}_1,...,{\chi}_{n-r})\otimes M$ over $k({\chi}_1,...,{\chi}_{n-r})$ in place of $M$ over $k$, it is thus equivalent to assume that $M$ is $n$-pure, that is $dim_k(M)<\infty$. In this case we have of course $M^{**}\simeq M$ and any associated prime ideal is maximum, that is $ass(M)\subset max(D)$. However, the reader must always keep in mind that the original module was pure and thus contained in its localization, that is to say no {\it simplification} is possible in the language of control theory.\\

As we have $soc(M)=\oplus soc_{\mathfrak{m}}(M)$ where the summation is {\it now} done on $ass(M)$ {\it only}, in order to dualize the short exact sequence $0\rightarrow soc(M)\rightarrow M\rightarrow M/soc(M)\rightarrow 0$, we need first dualize the various short exact sequences $0\rightarrow soc_{\mathfrak{m}}(M)\rightarrow M\rightarrow N\rightarrow 0$. However, if $S$ is a simple module, we have $\mathfrak{m}=ann(S)\subseteq ann(S^*)\subseteq ann(S^{**})=ann(S)$. Accordingly, the dual of a simple module isomorphic to $D/\mathfrak{m}$ is an isotypical simple module, because else it would have a proper factor module, the dual of which would be a proper submodule of $D/\mathfrak{m}$. We get therefore, again because of the injectivity of $D^*$, the short exact sequences $0\rightarrow N^*\rightarrow M^*\rightarrow top_{\mathfrak{m}}(M^*)\rightarrow 0$ and the desired dual sequence is finally obtained by introducing the intersection $rad(M^*)=\cap N^*$ for the various $\mathfrak{m}\in ass(M)$. A key result is provided by the following theorem which is not evident at all, even on very elementary examples, and provides a link between the socle of a module and the top of the corresponding system.\\

\noindent
{\bf THEOREM 5.2.2}: $N^*\simeq \mathfrak{m}M^*$ and the previous short exact sequence relative to $\mathfrak{m}$ is isomorphic to the short exact sequence $0\rightarrow \mathfrak{m}M^*\rightarrow M^*\rightarrow D/ \mathfrak{m}\otimes M^*\rightarrow 0$.\\
\noindent
{\it Proof}: As the second result is just obtained by tensoring with $M^*$ the short exact sequence $0\rightarrow \mathfrak{m}\rightarrow D\rightarrow D/ \mathfrak{m}\rightarrow 0$, it just remains to prove the first one.\\
For this, let us set $\mathfrak{m}=(a_1,...,a_t)=\sum Da$ and use the notations of the preceding lemma. If we introduce $I: a=J(a)\subset F$ and introduce the corresponding short exact sequence $0\rightarrow J(a)\rightarrow F\rightarrow N(a)\rightarrow 0$ for each generator $a\in \mathfrak{m}$, we have $aJ(a)\subseteq I$ by definition and the multiplication by $a$ thus induces a monomorphism $0\rightarrow N(a)\stackrel{a}{\rightarrow } M$. By duality, we have the epimorphism $M^*\stackrel{a}{\rightarrow }N(a)^*\rightarrow 0$ and obtain therefore $N(a)^*=aM^*$. Finally, if $I:\mathfrak{m}=J\subset F$, we have of course $J=\cap J(a)$ where the intersection is taken on the various generators of $\mathfrak{m}$ and an exact sequence $0\rightarrow J\rightarrow F\rightarrow \oplus N(a)$ inducing a monomorphism $0\rightarrow N\rightarrow \oplus N(a)$ if we define $N$ by the short exact sequence $0\rightarrow J\rightarrow F\rightarrow N\rightarrow 0$. Moreover, the inclusion $J\subseteq J(a), \forall a\in \mathfrak{m}$ induces an epimorphism $N\rightarrow N(a)\rightarrow 0$. Accordingly, we are {\it exactly} in the position to use Proposition 4.7 in order to get by duality the inclusion $0\rightarrow N(a)^*\rightarrow N^*$ and therefore $N^*=\sum N(a)^*=\sum aM^*=\mathfrak{m}M^*\subseteq M^*$. Finally, any nonzero element in the module on the right in the sequence of the theorem is killed by $\mathfrak{m}$ and admits a representative in $R$ which is not in $\mathfrak{m}R$. By definition of $hom$, it is the restriction of a section of $R$ to a simple submodule of $M$.\\ 
\hspace*{12cm}   Q.E.D.   \\

\noindent
{\bf COROLLARY 5.2.3}: We have the short exact sequence $0\rightarrow \cap \mathfrak{m} R\rightarrow R\rightarrow top(R)\rightarrow 0$ coming from the {\it chinese remainder theorem} by tensoring $D/ \cap {\mathfrak{m}}_i\simeq \oplus D/ {\mathfrak{m}}_i$ with $R$.\\

\noindent
{\bf COROLLARY 5.2.4}: $0\rightarrow R\rightarrow \oplus R_{\mathfrak{m}}$ projects onto $top(R)=\oplus top_{\mathfrak{m}}(R)=\oplus top(R)_{\mathfrak{m}}$.\\
\noindent
{\it Proof}: First of all, using the exactness of the localizing functor, we have $top(R)=\oplus D/\mathfrak{m}\otimes R$ with $D/\mathfrak{m}\otimes R\simeq R/\mathfrak{m}R\Rightarrow R_{\mathfrak{m}}/\mathfrak{m}R_{\mathfrak{m}}\simeq (R/\mathfrak{m}R)_{\mathfrak{m}}\simeq R/\mathfrak{m}R\simeq top_{\mathfrak{m}}(R)$. Indeed, $\mathfrak{m}\in max(D)\Rightarrow (\mathfrak{m},s)=D, \forall s\in D-\{\mathfrak{m}\}$ and $ \exists t\in D-\{\mathfrak{m}\}, a\in \mathfrak{m}$ with $st+a=1$. Accordingly, $\forall f\in R$, then $\frac{1}{s}f=\frac{st+a}{s}=tf+\frac{a}{s}f$ and we can therefore take out the denominators when localizing. Finally, as ${\mathfrak{m}}{'}+{\mathfrak{m}}{"}=D, \forall {\mathfrak{m}}{'},{\mathfrak{m}}{"}\in ass (M)$, we have similarly $(R/{\mathfrak{m}}{'}R)_{{\mathfrak{m}}{"}}=0$.\\
\hspace*{12cm}   Q.E.D.   \\

This corollary allows one to use Nakayama{'}s lemma in order to look for the generators of $R$ because $M$ and thus $R$ are finitely generated over $k$ and thus over $D$ by assumption. The following theorem, which is a straight consequence of ([12], IV, \S 2, p 104-109), constitutes the {\it secret} of Macaulay ([13], \S 82, end p 91) and explains the reason for introducing the (inverse) system.\\

\noindent
{\bf THEOREM 5.2.5}: When $M$ is $n$-pure, the minimum number of generators of $R$ is equal to the maximum number of isotypical components that can be found among the various components of the socle of $M$ or of the top of $R=M^*$, that is $max_{\mathfrak{m}\in ass(M)}\{dim _{D/\mathfrak{m}}soc_{\mathfrak{m}}(M)\}$.\\

\noindent
5.3)  As the examples in this subsection will clearly show out, it is important to notice that the number of generators is related to the localized module/system and {\it not} to the original module/system as we shall exhibit situations needing two generators even though $max(D)\cap ass(M)=\emptyset$. Therefore, in this last subsection, we shall explain the way followed by Macaulay ([13], \S 79, p 89) in order to get back informations on the original system from results on the localized one. For simplicity the index $k$ of the unknowns will not be written down.\\
Setting $\chi=({\chi}',{\chi}")$ with ${\chi}'=({\chi}_1,...,{\chi}_{n-r})$ and ${\chi}"=({\chi}_{n-r+1},...,{\chi}_n)$ while introducing similarly $\mu=(\mu{'},\mu{"})$ with ${\mu}'=({\mu}_1,...,{\mu}_{n-r})$ and $ {\mu}"=({\mu}_{n-r+1},...,{\mu}_n)$, we obtain the localized system by substituting $y_{\mu}=y_{(\mu{'},\mu{"})}={\chi}_{\mu{'}}y_{\mu{"}}$ in the original PD equations. However, if we start from an involutive system, the corresponding localized system may be still involutive with full classes but with quite different features, for example no longer homogeneous if the original system is homogeneous. In order to manage with a solved form, we have the following technical result found by Macaulay ([13], \S 78, p 88 (A) and \S 79, p 89 (B)).\\

\noindent
{\bf PROPOSITION 5.3.1}: The localized system is $k({\chi}')\otimes R $.\\
\noindent
{\it Proof}: As $M$ is finitely presented, there is the abstract isomorphism ([22], Th 3.84, p 107):\\
\[      hom_{k({\chi}')}(k({\chi}')\otimes M,k({\chi}'))\simeq k({\chi}' )\otimes hom_k(M,k)=k({\chi}')\otimes R   \]
where we recall that $M$ can be identified with its image in $k({\chi}')\otimes M$ as $M$ is r-pure. However, {\it in actual practice}, it is not evident at all to discover that {\it a single determinant in} $k({\chi}')$ allows to provide modular equations with coefficients polynomials in $k[{\chi}']$. In fact, all principal jets ($pri$) of order $\geq q$ and class $\geq n-r+1$ can be expressed from parametric jets ($par$) of the original system and can therefore be expressed by means of {\it finite} linear combinations of the jets of the localized system with order $\geq q-1$ with coefficients in $k[{\chi}']$ ({\it care}). However, these latter jets can themselves be linearly dependent through a finite number of equations of order $\leq q-1$ . Solving these finitely many equations with respect to principal jets of order $\leq q-1$ may therefore bring a determinant $c({\chi}')\in k[{\chi}']$. Accordingly, any modular equation of the localized system can be written in the form $E\equiv c({\chi}')a^{para}+\sum b({\chi}')a^{pri}=0$ with $b,c\in k[{\chi}']$ and we get therefore a finite number of modular equations of the form $E\equiv c_{{\mu}"}({\chi}')a^{{\mu}"}\equiv c^{{\lambda}'}_{{\mu}"}{\chi}_{{\lambda}'}a^{{\mu}"}=0$ called $r$-{\it dimensional modular equations} by Macaulay, with an inequality $\mid\lambda{'}\mid-\mid\mu{"}\mid\leq \delta$ for a certain relative integer $\delta$ and equality for homogeneous systems. Indeed, any PD equation is of the form $y^{pri}_{({\lambda}',{\lambda}")}\in \sum k y^{par}_{({\mu}',{\mu}")}$ with $\mid{\lambda}'\mid+\mid{\lambda}"\mid\geq \mid{\mu}'\mid+\mid{\mu}"\mid $. By localization, we get ${\chi}_{{\lambda}'}y^{pri}_{{\lambda}"}\in \sum c {\chi}_{{\mu}'}y^{par}_{{\mu}"}$. Setting $y^{par}_{{\mu}"}=1$ and the other parametric jets equal to zero, we obtain a modular equation of the form $E\equiv {\chi}_{{\lambda}'}a^{{\mu}"}_{par}+\sum c {\chi}_{{\mu}'}a^{{\lambda}"}_{pri}$ and thus $\mid{\mu}'\mid-\mid{\lambda}"\mid\leq \mid{\lambda}'\mid-\mid{\mu}"\mid\leq \delta$ as the number of parametric jets of order $\leq q-1$ in the localized system is finite. As no {\it simplification} may exist, that is $t_{r-1}(M)=0$ for the original module, one just needs to set $a^{\mu{"}}={\chi}_{\mu{'}}a^{(\mu{'},\mu{"})}$ in order to get $E\equiv c^{{\lambda}'}_{{\mu}"}{\chi}_{{\lambda}'+{\mu}'}a^{({\mu}',{\mu}")}=0$ and the so-called $n$-{\it dimensional modular equations} $E_{\alpha{'}}\equiv {\sum}_{\lambda{'}+\mu{'}=\alpha{'}}c^{\lambda{'}}_{\mu{"}}a^{(\mu{'},\mu{"})}=0$. \\
\hspace*{12cm}  Q.E.D.   \\

\noindent
{\bf EXAMPLE 5.3.2}:  {\it Purity is essential in the process of localization/delocalization}. In order to prove it, let us consider the very simple Example 2.3 of codimension 1 but not pure. Localizing, we get from the second PD equation ${\chi}_1y_2=0$ with a {\it simplification} leading to the new PD equation $y_2=0$ (We let the reader compare with the pure situation of Examples 2.2 and 3.15).\\

The following theorem on the way to generate the modular equations is the key result obtained by Macaulay ([13], \S 82, end of p 91).\\

\noindent
{\bf THEOREM 5.3.3}: There is a finite number of $r$-dimensional modular equations, a smaller number 
of $r$-dimensional modular equations of which all the others are derivates and an equal or still smaller number of $n$-dimensional modular equations of which all the others of an arbitrary order $q$ are derivates.\\
\noindent
{\it Proof}: In $E_{\alpha{'}}$ we have $\mid\mu\mid=\mid\mu{'}\mid + \mid\mu{"}\mid=\mid\alpha{'}\mid-\mid\lambda{'}\mid + \mid\mu{"}\mid\geq \mid\alpha{'}\mid-\delta$. As for the derivates, we have $d_{\gamma{"}}E\equiv c^{\lambda{'}}_{\mu{"}}{\chi}_{\lambda{'}+\mu{'}}a^{(\mu{'},\mu{"}-\gamma{"})}=0$ with $a^{\mu-1_i}=0$ if ${\mu}_i=0$ or $a^{({\mu}_1,...,{\mu}_i-1,...,{\mu}_n)}$ if ${\mu}_i\geq 1$ ([13], \S 60, p 69). \\
We have therefore $(d_{\gamma{"}}E)_{\beta{'}}\equiv d_{\gamma{"}}(E_{\beta{'}})\equiv d_{\gamma{"}}E_{\beta{'}}\equiv {\sum}_{\lambda{'}+\mu{'}=\beta{'}} c^{\lambda{'}}_{\mu{"}}a^{(\mu{'},\mu{"}-\gamma{"})}=0$ with $\mid\gamma{"}\mid \leq \tau$ since there are only a finite number of linearly independent derivates of the $r$-dimensional equations as the localized system is a finite dimensional differential vector space over $k({\chi}')$. \\
More precisely, any $\gamma{"}$-derivates is of the form $d_{{\gamma}"}E_{\beta{'}}=0 $ with $\lambda{'}\leq \beta{'},\gamma{"}\leq \mu{"}, \mid{\gamma}"\mid\leq \tau$ and where ${\beta}',{\gamma}"$ are {\it fixed} multi-indices. \\
Let us finally consider all the modular equations of order $q$ that can be obtained as derivates, that is all the $d_{{\gamma}"}E_{\beta{'}}$ with $q\geq \mid\beta{'}\mid-\mid\lambda{'}\mid+\mid\mu{"}\mid-\mid\gamma{"}\mid\geq \mid\beta{'}\mid -\tau -\delta$ that is to say ${\beta{'}}_i\leq \mid\beta{'}\mid \leq q+\delta+\tau, \forall i=1,...,n-r$. Accordingly, every modular equation of order $q$ is a derivate of a certain $E_{\alpha{'}}$ for a fixed $\alpha{'}$ if ${\alpha{'}}_i\geq q+\delta+\tau=q', \forall i=1,...,n-r$.\\
Let us explain this fact in the simple 2-dimensional situation ${\beta}'=({\beta}^{'}_1\geq 0,{\beta}^{'}_2\geq 0)$ arising when $r=n-2$. Any modular equation $E_{{\beta}'}$ provides a point ${\beta}'$ in this quadrangle and all the modular equations of order $q$ come therefore from points contained in the triangle made by the two axes and the straight line ${\beta}^{'}_1+{\beta}^{'}_2=q'$ which can all be obtained as derivates of $E_{(q',q')}$.\\
\hspace*{12cm}  Q.E.D.   \\

\noindent
{\bf EXAMPLE 5.3.4}: Coming back to Example 5.1.2, we have $par=\{y,y_1,y_2,y_{11}\}\Rightarrow M\simeq ky+ky_1+ky_2+ky_{11}$ and thus $f=(1,0,0,0)\rightarrow E_1\equiv a^0=0, f=(0,1,0,0)\rightarrow E_2\equiv a^1=0, f=(0,0,1,0)\rightarrow E_3\equiv a^2=0, f=(0,0,0,1)\rightarrow E_4\equiv a^{11}=0$. We have $soc(M)\simeq ky_2+ ky_{11}$ with two isotypical components both killed by $\mathfrak{m}=(d_1,d_2)$ and  $top(R)=\{E_3,E_4\}$ provides two generators for the $4$-dimensional differential vector space $R$ as we have indeed $d_2E_3=E_1, d_1E_4=E_2$, that is to say $\mathfrak{m}R$ is generated by $\{E_1,E_2\}$ in agrement with Nakayama's lemma (compare to [6], p 526).\\

\noindent
{\bf EXAMPLE 5.3.5}: Coming back to Example 5.1.3, we have $par=\{y, y_1, y_2, y_3, y_{11}, y_{12}, y_{13}, y_{111}\}$. We have $soc(M)\simeq ky_{111}$ killed by $\mathfrak{m}=(d_1,d_2,d_3)$ and $top(R)=\{E\}$ with $E\equiv a^{111}+a^{123}=0$ provides a unique generator for the $8$-dimensional differential vector space $R$ as we have $d_1E\equiv a^{11}+a^{23}=0, d_2E\equiv a^{13}=0, d_3E\equiv a^{12}=0, ..., d_{111}E\equiv a^0=0$ and a way to generate $\mathfrak{m}R$. It is {\it remarkable} that $8=2^3=2^n$ is a general combinatorial result proved by Macaulay ([13], \S 58, p 79, \S 84, p 92). A similar simpler situation is met with $n=2$ and $y_{22}=0, y_{12}-y_{11}=0$, leading to $E\equiv a^{11}+a^{12}=0$ or with $n=3$ and $y_{33}-y_{11}=0, y_{23}=0, y_{22}-y_{11}=0$ leading to $E\equiv a^{111}+a^{122}+a^{133}=0$ (compare to [13], p 81).\\

\noindent
{\bf EXAMPLE 5.3.6}: Coming back to Example 3.14 which needs a partial localization with $k({\chi}')=k({\chi}_1,{\chi}_2,{\chi}_3)$, the localized system is $y^1_4=0, y^2_4=0, y^3_4=0, {\chi}_3y^3+{\chi}_2y^2+{\chi}_1y^1=0$. Clearly $k({\chi}')\otimes M\simeq k({\chi}')y^1+k({\chi}')y^2$ is a semi-simple module with two isotypical components both killed by $\mathfrak{m}=(d_4)$. Accordingly, $(1,0)\rightarrow E_1\equiv {\chi}_3a^0_1-{\chi}_1a^0_3=0, (0,1)\rightarrow E_2\equiv {\chi}_3a^0_2-{\chi}_2a^0_3 =0$ provides the two generators of the localized system, even though $max(D)\cap ass(M)=\emptyset$ in this case. We notice that the determinant $c({\chi}')={\chi}_3$ is unavoidable. Delocalizing, we get ${\chi}_1{\chi}_3\rightarrow a^1_1-a^3_3=0, {\chi}_2{\chi}_3\rightarrow a^2_2-a^3_3=0$ and so on, in agrement with the general theory for the original system.\\

\noindent
{\bf EXAMPLE 5.3.7}: With $n=3, m=1, q=2, k=\mathbb{Q}$, the module defined by the homogeneous involutive system $y_{33}=0, y_{23}-y_{13}=0, y_{22}-y_{12}=0$ is 2-pure. Setting $k({\chi}')=k({\chi}_1)$, the corresponding localized system  $y_{33}=0, y_{23}-{\chi}_1y_3=0, y_{22}-{\chi}_1y_2=0$ is again involutive with $par=\{y, y_2, y_3\}$. We obain therefore $(1,0,0)\rightarrow E_1\equiv a^0=0, (0,1,0)\rightarrow E_2\equiv a^2+{\chi}_1a^{22}+({\chi}_1)^2a^{222}+... =0, (0,0,1)\rightarrow E_3\equiv a^3+{\chi}_1a^{23}+({\chi}_1)^2a^{223}+ ... =0$. We notice that ${\mathfrak{m}}_1=(d_3,d_2-{\chi}_1)$ kills $y_3$ while ${\mathfrak{m}}_2=(d_3,d_2)$ kills $y_2-{\chi}_1y$, each maximum ideal in $k({\chi}_1)[d_2,d_3]$ leading to a unique isotypical component. Denoting simply by $M$ the localized module and by $R$ the corresponding system, we have the short exact sequence $0\rightarrow soc_{{\mathfrak{m}}_2}(M)\rightarrow M\rightarrow N_2\rightarrow 0$ (care to the notation) and the dualizing short exact sequence $0\rightarrow N^*_2\rightarrow R\rightarrow top_{{\mathfrak{m}}_2}(R)\rightarrow 0$. Here, $N_2$ is obtained by adding $y_2-{\chi}_1y=0$ to the equations of the localized system and we obtain the subsystem $N^*_2=f_1(E_1+{\chi}_1E_2)+f_3E_3\subset R=f_1E_1+f_2E_2+f_3E_3$. We check the relations:\\
\[ d_2E_1=0, d_3E_1=0, d_2E_2=E_1+{\chi}_1E_2, d_3E_2=0, d_2E_3={\chi}_1E_3, d_3E_3=E_1+{\chi}_1E_2\]
transforming $R$ into a 3-dimensional differential vector space over $k({\chi}_1)$ and obtain therefore $d_2R=f_2(E_1+{\chi}_1E_2)+{\chi}_1f_3E_3, d_3R=f_3(E_1+{\chi}_1E_2)$, that is we check directly $N^*_2={\mathfrak{m}}_2R$ and could check similarly $N^*_1={\mathfrak{m}}_1R$, a result highly not evident at first sight. According to the general theory, there should be one generator only and we may choose $E=E_2+E_3$ in order to generate $R$ as we have indeed the three linearly independent relations:\\
\[   E=E_2+E_3, d_2E-{\chi}_1E=E_1, d_3E=E_1+{\chi}_1E_2\]
allowing to determine $E_1,E_2,E_3$ from the derivates of $E$. The system being homogeneous, we have $q'=q+\delta+\tau=2-1+1=2$. As we have $E\equiv 1(a^2+a^3)+ {\chi}_1(a^{12}+a^{13}+a^{22}+a^{23}) +({\chi}_1)^2(a^{112}+a^{122}+a^{222}+a^{113}+a^{123}+a^{223})+ ... =0$, it is easy to check that the single modular equation $E_{11}\equiv 
a^{112}+a^{122}+a^{222}+a^{113}+a^{123}+a^{223}=0$ generates $a^{11}=0, a^{13}+a^{23}=0, a^{12}+a^{22}=0, a^1=0, a^2=0, a^3=0, a^0=0$ successively. Hence all the modular equations up to order 2 are generated by a single modular equation at order 3, a result not evident at first sight.\\

\noindent
{\bf EXAMPLE 5.3.8}: Looking back to Example 2.3, the primary decomposition brings the two subsystems $R{'}=\{y_2=0\rightarrow a^0=0,a^1=0,a^{11}=0,...\}$ and $R{"}=\{y_{22}=0,y_{12}=0,y_{11}=0\rightarrow a^0=0,a^1=0,a^2=0\}$ with $R=R{'}+R{"}$. One clearly needs {\it two} generators in order to generate any $R_q$, say $\{a^2=0,a^{111}=0\}$ for $q=3$. More generally, the involutive system $y_{22}=0, y_{12}-ay_2=0, y_{11}-ay_1=0$ depending on the constant parameter $a$ needs 1 generator if $a\neq 0$ (generic case just studied) but 2 if $a=0$. The situation is similar with the system $y^1_{xx}-ay^1=0, y^2_x=0$ presented at the end of the introduction when $a=0$ and $a=1$ ([13], \S 83, p 91). Such a result proves that the identifyability of a system may depend on the parameters involved and refines the classification of systems or modules presented in [20].\\

\noindent
{\bf EXAMPLE 5.3.9}: Coming back to the Example 2.2 of Macaulay, the inhomogeneous involutive system $y_{33}=0, y_{23}=0, y_{22}=0, y_{13}-y_2=0$ has the unique generating 3-dimensional modular equation $E\equiv 1a^3+{\chi}_1(a^2+a^{13})+... =0$ in a coherent way with ([13], $\S$72).\\

\noindent
{\bf EXAMPLE 5.3.10}: With now $n=4$, let us consider the 2-pure module defined by the homogeneous involutive system $y_{44}=0, y_{34}=0, y_{33}=0, y_{24}-y_{13}=0$. Using ${\chi}'=({\chi}_1,{\chi}_2)$, we get the single generating 2-dimensional equation $E\equiv {\chi}_1a^4+{\chi}_2a^3=0$ and the corresponding 4-dimensional modular equation $E\equiv {\chi}_1a^4+{\chi}_2a^3+({\chi}_1)^2a^{14}+{\chi}_1{\chi}_2(a^{24}+a^{13})+({\chi}_2)^2a^{23}+ ...=0$ with $d_4E_{112}\equiv a^{12}=0 $ for example and $q'=q+1$ though $q'=q$ is sufficient here. The module $L$ with projective dimension 2 can be defined by the involutive system $z_4=0, y_4-z_1=0, z_3=0, y_3-z_2=0$ for $(y,z)$ and there is a strict inclusion $M\subset L$ obtained by eliminating $z$.\\

\noindent
6)  {\bf CONCLUSION}:\\
We summarize the way leading to revisit the {\it inverse system} of Macaulay by using modern techniques of {\it algebraic analysis}, namely differential geometric arguments for studying the system instead of the module. \\
The main purpose is to find generators for the {\it differential system} dual to the {\it differential module}. For this, the only way known in the literature is to control the generators of the system from the generators of its {\it top} by using Nakayama's lemma. Again by {\it duality}, this amounts to count the isotypical components of the {\it socle} of the module. Meanwhile, the key idea is to decompose the system into subsystems instead of using a primary decomposition of the module in order to deal with {\it pure} modules, a concept generalizing the {\it unmixedness assumption} of Macaulay. However, the original system is {\it not} in general finitely generated and it is therefore {\it essential} to use a {\it partial localization} in order to deal with a finite dimensional localized system.\\
The present approach avoids the abstract systematic use of the {\it injective hull} by Oberst and opens a new way towards effective {\it computer algebra} packages for studying {\it identifiability} in engineering sciences.\\
It is thus remarkable that Macaulay had the intuition of these techniques as early as in 1916 and we express our deep gratitude to his work.\\

\noindent
{\bf BIBLIOGRAPHY}\\
 
\noindent
[1] I. ASSEM, Alg\`{e}bres et Modules, Masson, Paris, 1997.\\
\noindent
[2] J.E. BJORK, Analytic D-modules and Applications, Kluwer, 1993.\\
\noindent
[3] BOURBAKI, Alg\`{e}bre Commutative, Chapitre 1 \`{a} 4, Masson, Paris, 1985.\\
\noindent
[4] BOURBAKI, Alg\`{e}bre, Chapitre 10, Alg\`{e}bre commutative, Masson, Paris, 1980.\\
\noindent
[5] B. BUCHBERGER, Ein Algorithmus zum Auffinden der Basiselemente des Restklassenringes nach einem Multidimensionalen Polynomideal, PhD thesis (thesis advisor W. Gr\"obner), University of Innsbruck, Austria, 1965. English translation: An Algorithm for Finding the Basis Elements in the 
Residue Class Ring Modulo a Zero Dimensional Polynomial Ideal, Journal of Symbolic Computations, Special Issue on logic, Mathematics and Computer Sciences: Interactions, Vol 14, Nb 34, 2006, 475-511. \\ 
\noindent
[6] D. EISENBUD, Commutative Algebra With a View Towards Algebraic Geometry, Graduate Texts in Math 150, Springer, 1996.\\
\noindent
[7] V.P. GERDT, Y.A. BLINKOV, Minimum Involutive Bases, Mathematics and
Computers in Simulations, 45, 1998, 543-560.\\
\noindent
[8] W. GR\"{O}BNER, \"{U}ber die Algebraischen Eigenschaften der Integrale 
von Linearen Differentialgleichungen mit Konstanten Koeffizienten, 
Monatsh. der Math., 47, 1939, 247-284.\\
\noindent
[9] M. JANET, Sur les Syst\`emes aux d\'eriv\'ees partielles, 
Journal de Math., 8, 3, 1920, 65-151.\\
\noindent
[10] E.R. KALMAN, Y.C. YO, K.S. NARENDA, Controllability of Linear Dynamical
Systems, Contrib. Diff. Equations, 1, 2, 1963, 189-213.\\
\noindent
[11] M. KASHIWARA, Algebraic Study of Systems of Partial Differential 
Equations, M\'emoires de la Soci\'et\'e Math\'ematique de France 63, 1995, 
(Transl. from Japanese of his 1970 Master's Thesis).\\
\noindent
[12] E. KUNZ, Introduction to Commutative Algebra and Algebraic Geometry, 
Birkh\"{a}user, 1985.\\
\noindent
[13] F. S. MACAULAY, The Algebraic Theory of Modular Systems, Cambridge Tracts 19, Cambridge University Press, London, 1916; Reprinted by Stechert-Hafner Service Agency, New York, 1964.\\
\noindent
[14] D.G. NORTHCOTT: Lessons on Rings, Modules and Multiplicities, Cambridge University Press, 1968.\\
\noindent
[15] U. OBERST, Multidimensional Constant Linear Systems, Acta Appl. Math., 20, 1990, 1-175.\\
\noindent
[16] V.P. PALAMODOV, Linear Differential Operators with Constant Coefficients,
Grundlehren der Mathematischen Wissenschaften 168, Springer, 1970.\\
\noindent
[17] J.-F. POMMARET, Partial Differential Equations and Group Theory,New Perspectives for Applications, Mathematics and its Applications 293, Kluwer, 1994.\\
\noindent
[18] J.-F. POMMARET, Partial Differential Control Theory, Kluwer, 2001, 957 pp.\\
(http://cermics.enpc.fr/$\sim$pommaret/home.html)\\
\noindent
[19] J.-F. POMMARET, Algebraic Analysis of Control Systems Defined by Partial Differential Equations, in Advanced Topics in Control Systems Theory, Lecture Notes in Control and Information Sciences 311, Chapter 5, Springer, 2005, 155-223.\\
\noindent
[20] J.-F. POMMARET, Gr\"{o}bner Bases in Algebraic Analysis: New perspectives for applications, Radon Series Comp. Appl. Math 2, 1-21, de Gruyter, 2007.\\
\noindent
[21] A. QUADRAT: http://wwwb.math.rwth-aachen.de/OreModules  \\
http://www.risc.uni-linz.ac.at/about/conferences/aaca09/   , in particular ... /ModuleTheoryI.pdf   and   ... /ModuleTheoryII.pdf  \\
\noindent
[22] J.J. ROTMAN, An Introduction to Homological Algebra, Pure and Applied 
Mathematics, Academic Press, 1979.\\
\noindent
[23] W. M. SEILER, Involution: The Formal Theory of Differential Equations and its Applications to Computer Algebra, Springer, 2009, 660 pp.\\
\noindent
[24] D.C. SPENCER, Overdetermined Systems of Partial Differential Equations, 
Bull. Amer. Math. Soc., 75, 1965, 1-114.\\

\end{document}